\documentclass[twoside]{article}
\usepackage{txmac,a4,
amssymb,pifont,%
case,%
nocaphead,%
epsf,%
myrot,%
mypic,times,mathptm}

\advance\oddsidemargin by -1.9cm
\advance\evensidemargin by -1.9cm
\advance\textwidth by 3.8cm

\def\mynewtheo#1#2{%
\newtheorem{@#1}{#2}[section]%
\newenvironment{#1}{\begin{@#1}\rm}{\end{@#1}}}

\makeatletter

\newtheorem{@clm}{Claim}

\def\the@clm{\arabic{@clm}.}

\mynewtheo{lemma}{Lemma}
\mynewtheo{slemma}{Sublemma}
\mynewtheo{exer}{Exercise}
\mynewtheo{theo}{Theorem}
\mynewtheo{rem}{Remark}
\mynewtheo{defi}{Definition}
\mynewtheo{conj}{Conjecture}
\mynewtheo{corr}{Corollary}
\mynewtheo{prop}{Proposition}
\mynewtheo{question}{Question}
\mynewtheo{exam}{Example}

\newenvironment{theorem}{\begin{theo}}{\end{theo}}

\begin{document}

\newenvironment{eqn}{\begin{equation}}{\end{equation}\@ignoretrue}

\newenvironment{myeqn*}[1]{\begingroup\def\@eqnnum{\reset@font\rm#1}%
\xdef\@tempk{\arabic{equation}}\begin{equation}\edef\@currentlabel{#1}}
{\end{equation}\endgroup\setcounter{equation}{\@tempk}\ignorespaces}

\newenvironment{myeqn}[1]{\begingroup\let\eq@num\@eqnnum
\def\@eqnnum{\bgroup\let\r@fn\normalcolor 
\def\normalcolor####1(####2){\r@fn####1#1}%
\eq@num\egroup}%
\xdef\@tempk{\arabic{equation}}\begin{equation}\edef\@currentlabel{#1}}
{\end{equation}\endgroup\setcounter{equation}{\@tempk}\ignorespaces}

\newenvironment{myeqn**}{\begin{myeqn}{(
\theequation)\es\es\mbox{\qed}}\edef\@currentlabel{\theequation}}
{\end{myeqn}\stepcounter{equation}}

\newcommand{\mybin}[2]{\text{$\Bigl(\begin{array}{@{}c@{}}#1\\#2%
\end{array}\Bigr)$}}
\newcommand{\mybinn}[2]{\text{$\biggl(\begin{array}{@{}c@{}}%
#1\\#2\end{array}\biggr)$}}

\def\overtwo#1{\mbox{\small$\mybin{#1}{2}$}}
\newcommand{\mybr}[2]{\text{$\Bigl\lfloor\mbox{%
\small$\displaystyle\frac{#1}{#2}$}\Bigr\rfloor$}}
\def\mybrtwo#1{\mbox{\mybr{#1}{2}}}

\def\myfrac#1#2{\raisebox{0.2em}{\small$#1$}\!/\!\raisebox{-0.2em}{\small$#2$}}

\author{A. Stoimenow\footnotemark[1]\\[2mm]
\small Research Institute for Mathematical Sciences, \\
\small Kyoto University, Kyoto 606-8502, Japan\\
\small e-mail: {\tt stoimeno@kurims.kyoto-u.ac.jp}\\
\small WWW: {\hbox{\web|http://www.kurims.kyoto-u.ac.jp/~stoimeno/|}}
}

{\def\thefootnote{}
\footnotetext[1]{This is a preprint. I would be grateful
  for any comments and corrections. Current
  version: \today\ \ \ First version: \makedate{20}{8}{2004}}
\def\thefootnote{\fnsymbol{footnote}}
\footnotetext[1]{Financial support by the 21st Century COE Program
  is acknowledged.}
}

\title{%
\large\bf \uppercase{Alexander polynomials and}\\[2mm]
\uppercase{hyperbolic volume of arborescent links}%
}

\date{}

\maketitle

\let\vn\varnothing
\let\ay\asymp
\let\pa\partial
\let\ap\alpha
\let\be\beta
\let\bt\beta
\let\Dl\Delta
\let\nb\nabla
\let\Gm\Gamma
\let\gm\gamma
\let\de\delta
\let\dl\delta
\let\eps\epsilon
\let\lm\lambda
\let\Lm\Lambda
\let\sg\sigma
\let\vp\varphi
\let\om\omega
\let\diagram\diag
\let\ol\overline
\let\sm\setminus
\let\tl\tilde
\def\ncap{\not\mathrel{\cap}}
\def\sgn{\text{\rm sgn}\,}
\def\spn{\mathop {\operator@font span}}
\def\et{\mathop {\operator@font ext}}
\def\ir{\mathop {\operator@font int}}
\def\Md{\max\deg}
\def\md{\min\deg}
\def\mc{\max\cf}
\def\vol{\text{\rm vol}\,}
\def\hra{\hookrightarrow}
\def\Lra{\Longrightarrow}
\def\lra{\longrightarrow}
\def\so{\Rightarrow}
\def\So{\Longrightarrow}
\def\nin{\not\in}
\let\ds\displaystyle
\let\llra\longleftrightarrow
\let\tg\triangle
\let\bd\partial
\let\reference\ref
\let\wh\widehat
\let\wt\widetilde
\let\op\oplus
\def\TM{$^\text{\raisebox{-0.2em}{${}^\text{TM}$}}$}
\def\lb{\linebreak[0]}
\def\lz{\lb\verb}
\def\ssim{\stackrel{\ds
\sim}{\vbox{\vskip-0.2em\hbox{$\scriptstyle *$}}}}
\def\btp{\bar t_2'}
\let\es\enspace

\long\def\@makecaption#1#2{%
   \vskip \abovecaptionskip 
   {\let\label\@gobble
   \let\ignorespaces\@empty
   \xdef\@tempt{#2}%
   }%
   \ea\@ifempty\ea{\@tempt}{%
   \sbox\@tempboxa{%
      \fignr#1#2}%
      }{%
   \sbox\@tempboxa{%
      {\fignr#1:}\capt\ #2}%
      }%
   \ifdim \wd\@tempboxa >\captionwidth {%
      \rightskip=\@captionmargin\leftskip=\@captionmargin
      \unhbox\@tempboxa\par
     }%
   \else
      \centerline{\box \@tempboxa}%
   \fi
   \vskip \belowcaptionskip
   }%
\def\fignr{\small\sffamily\bfseries}%
\def\capt{\small\sffamily}%


\newdimen\@captionmargin\@captionmargin2cm\relax
\newdimen\captionwidth\captionwidth0.8\hsize\relax

\def\rottab#1#2{%
\expandafter\advance\csname c@table\endcsname by -1\relax
\centerline{%
\rbox{\hsize=\vsize\relax\centerline{\vbox{\setbox1=\hbox{#1}%
\centerline{\mbox{\hbox to \wd1{\hfill\mbox{\vbox{{%
\caption{#2}}}}\hfill}}}%
\vskip9mm
\centerline{
\mbox{\copy1}}}}%
}%
}%
}

\newcounter{local}
\newenvironment{mylist}[1]{
\bgroup
\def\myitem{\item[\stepcounter{local}#1{local})\hfil]%
\xdef\@tempa{\number\c@local}
}
\ea\gdef\ea\@currentlabel\ea{#1{local}}
\begin{list}{}{\setcounter{local}{0}%
\itemsep1mm\relax 
\leftmargin9mm\relax
\topsep1mm\relax \labelwidth5mm }}{\end{list}\egroup}

\def\eqref#1{(\protect\ref{#1})}

\def\proof{\@ifnextchar[{\@proof}{\@proof[\unskip]}}
\def\@proof[#1]{\noindent{\bf Proof #1.}\es}

\def\hint{\noindent Hint: }
\def\problem{\noindent{\bf Problem.} }

\def\@mt#1{\ifmmode#1\else$#1$\fi}
\def\qed{\hfill\@mt{\Box}}
\def\qqed{\hfill\@mt{\Box\es\Box}}

\def\@curvepath#1#2#3{%
  \@ifempty{#2}{\piccurveto{#1 }{@stc}{@std}#3}%
    {\piccurveto{#1 }{#2 }{#2 #3 0.5 conv}
    \@curvepath{#3}}%
}
\def\curvepath#1#2#3{%
  \piccurve{#1 }{#2 }{#2 }{#2 #3 0.5 conv}%
  \picPSgraphics{/@stc [ #1 #2 -1 conv ] $ D /@std [ #1 ] $ D }%
  \@curvepath{#3}%
}

\def\@opencurvepath#1#2#3{%
  \@ifempty{#3}{\piccurveto{#1 }{#1 }{#2 }}%
    {\piccurveto{#1 }{#2 }{#2 #3 0.5 conv}\@opencurvepath{#3}}%
}
\def\opencurvepath#1#2#3{%
  \piccurve{#1 }{#2 }{#2 }{#2 #3 0.5 conv}%
  \@opencurvepath{#3}%
}

\def\cA{{\cal A}}
\def\cU{{\cal U}}
\def\cC{{\cal C}}
\def\cP{{\cal P}}
\def\fg{{\frak g}}
\def\Kr{\mathop{\operator@font Kr}}
\def\diam{\mathop{\operator@font diam}}
\def\cZ{{\cal Z}}
\def\cD{{\cal D}}
\def\bQ{{\Bbb Q}}
\def\bR{{\Bbb R}}
\def\cE{{\cal E}}
\def\bZ{{\Bbb Z}}
\def\bN{{\Bbb N}}

\def\rato#1{\hbox to #1{\rightarrowfill}}
\def\hookrato#1{\hbox to
#1{$\lhook\joinrel$\rightarrowfill}}
\def\lhra{\hookrato{1.8em}}
 
\def\namedarrow#1{{\es
\setbox7=\hbox{F}\setbox6=\hbox{%
\setbox0=\hbox{\footnotesize
$#1$}\setbox1=\hbox{$\to$}%
\dimen@\wd0\advance\dimen@ by 0.66\wd1\relax
$\stackrel{\rato{\dimen@}}{\copy0}$}%
\ifdim\ht6>\ht7\dimen@\ht7\advance\dimen@ by
-\ht6\else
\dimen@\z@\fi\raise\dimen@\box6\es}}

\def\geni{\diag{3mm}{5}{1.7}{
  \piccirclearc{0.5 0}{0.5}{0 180}
  \piccirclearc{2.5 0}{0.5}{0 180}
  \piccirclearc{4.5 0}{0.5}{0 180}
}}

\def\genii{\diag{3mm}{5}{1.7}{
  \piccirclearc{1.5 0}{0.5}{0 180}
  \picellipsearc{1.5 0}{1.5 1.2}{0 180}
  \piccirclearc{4.5 0}{0.5}{0 180}
}}

\def\geniii{\diag{3mm}{5}{1.7}{
  \piccirclearc{0.5 0}{0.5}{0 180}
  \picellipsearc{3.5 0}{1.5 1.2}{0 180}
  \piccirclearc{3.5 0}{0.5}{0 180}
}}

\def\geniv{\diag{3mm}{5}{2.5}{
  \piccirclearc{1.5 0}{0.5}{0 180}
  \piccirclearc{3.5 0}{0.5}{0 180}
  \picellipsearc{2.5 0}{2.5 2}{0 180}
}}

\def\genv{\diag{3mm}{5}{2.5}{
  \piccirclearc{2.5 0}{0.5}{0 180}
  \picellipsearc{2.5 0}{2.5 2}{0 180}
  \picellipsearc{2.5 0}{1.5 1.2}{0 180}
}}

\def\twoegdr{\,{\diag{1.7em}{1}{0.8}{\picfillgraycol{0}
	\pictranslate{0 0.19}{
        \picfilledcircle{0 0.25}{0.08}{}
        \picfilledcircle{1 0.25}{0.08}{}
        \piccurve{0 0.25}{0.1 0.8}{0.9 0.8}{1 0.25}
        \piccurve{0 0.25}{0.1 -0.3}{0.9 -0.3}{1 0.25}
	}}}\,
}
\def\tmean{\,\diag{1.5em}{1}{0.8}{
   \piccircle{0.4 0.5 x}{0.37}{}
   \piclinedash{0.1}{0.05}
   \picline{0.4 0 x}{0.4 1 x}
   }\,}

\def\vrt#1{{\picfillgraycol{0}\picfilledcircle{#1}{0.06}{}}}
\def\vrtx#1#2{{\picfillgraycol{0}\picfilledcircle{#1}{#2}{}}}

\def\tycl#1#2#3#4{\vrt{#1}\vrt{#2}\vrt{#3}\vrt{#4}
\picline{#1}{#2}\picline{#2}{#3}\picline{#3}{#1}%
\picline{#2}{#4}\picline{#4}{#3}%
}

\def\trig#1#2#3{\vrtx{#1}{0.06}\vrtx{#2}{0.06}\vrtx{#3}{0.06}
\picline{#1}{#2}\picline{#2}{#3}\picline{#3}{#1}%
}

\def\ttycl#1#2#3#4#5{\vrtx{#1}{0.07}\vrtx{#2}{0.07}
\vrtx{#3}{0.07}\vrtx{#4}{0.07}\vrtx{#5}{0.07}
\picline{#1}{#2}\picline{#2}{#3}\picline{#3}{#1}%
\picline{#2}{#4}
\picline{#4}{#5}\picline{#5}{#3}%
}

\def\cycl#1#2#3#4{\vrt{#1}\vrt{#2}\vrt{#3}\vrt{#4}%
\picline{#1}{#2}\picline{#2}{#3}\picline{#3}{#4}\picline{#4}{#1}%
}

\def\xcycl#1#2#3#4#5#6{\vrt{#1}\vrt{#2}\vrt{#3}\vrt{#4}\vrt{#5}
\vrt{#6}\picline{#1}{#2}\picline{#2}{#3}\picline{#3}{#4}
\picline{#4}{#5}\picline{#5}{#6}\picline{#6}{#1}}

\def\llpoint#1{{\picfillgraycol{0}\picfilledcircle{#1}{0.15}{}}}
\def\lpoint#1{{\picfillgraycol{0}\picfilledcircle{#1}{0.08}{}}}
\def\point#1{{\picfillgraycol{0}\picfilledcircle{#1}{0.04}{}}}

\def\bysame{\same[\kern2cm]\,}

\def\br#1{\left\lfloor#1\right\rfloor}
\def\BR#1{\left\lceil#1\right\rceil}

\def\abstractname{}


\renewcommand{\section}{%
   \@startsection
         {section}{1}{\z@}{-1.5ex \@plus -1ex \@minus
-.2ex}%
               {1ex \@plus.2ex}{\large\bf}%
}
\renewcommand{\@seccntformat}[1]{\csname
the#1\endcsname .
\quad}

\def\bC{{\Bbb C}}
\def\bP{{\Bbb P}}

\def\epsfs#1#2{{\catcode`\_=11\relax\ifautoepsf\unitxsize#1\relax\else
\epsfxsize#1\relax\fi\epsffile{#2.eps}}}
\def\epsfsv#1#2{{\vcbox{\epsfs{#1}{#2}}}}
\def\vcbox#1{\setbox\@tempboxa=\hbox{#1}\parbox{\wd\@tempboxa}{\box
  \@tempboxa}}
\def\p{\epsfsv{4cm}}
\def\myss#1{\vcbox{#1}}

\def\eepsfs{\@ifnextchar[{\@eepsfs}{\@@eepsfs}}
\def\@eepsfs[#1]#2{\uu{\ifautoepsf\unitxsize#1\relax\else
\epsfxsize#1\relax\fi\epsffile{#2.eps}}}
\def\@@eepsfs#1{\uu{\epsffile{#1.eps}}}

\def\uu#1{\setbox\@tempboxa=\hbox{#1}\@tempdima=-0.5\ht\@tempboxa
\advance\@tempdima by 0.5em\raise\@tempdima\box\@tempboxa}

\def\@test#1#2#3#4{%
  \let\@tempa\go@
 
\@tempdima#1\relax\@tempdimb#3\@tempdima\relax\@tempdima#4\unitxsize\relax
  \ifdim \@tempdimb>\z@\relax
    \ifdim \@tempdimb<#2%
      \def\@tempa{\@test{#1}{#2}}%
    \fi
  \fi
  \@tempa
}

\def\go@#1\@end{}
\newdimen\unitxsize
\newif\ifautoepsf\autoepsftrue

\unitxsize4cm\relax
\def\epsfsize#1#2{\epsfxsize\relax\ifautoepsf
  {\@test{#1}{#2}{0.1 }{4   }
		{0.2 }{3   }
		{0.3 }{2   }
		{0.4 }{1.7 }
		{0.5 }{1.5 }
		{0.6 }{1.4 }
		{0.7 }{1.3 }
		{0.8 }{1.2 }
		{0.9 }{1.1 }
		{1.1 }{1.  }
		{1.2 }{0.9 }
		{1.4 }{0.8 }
		{1.6 }{0.75}
		{2.  }{0.7 }
		{2.25}{0.6 }
		{3   }{0.55}
		{5   }{0.5 }
		{10  }{0.33}
		{-1  }{0.25}\@end
		\ea}\ea\epsfxsize\the\@tempdima\relax
		\fi
		}

\let\lra\longrightarrow
\let\sm\setminus
\let\eps\varepsilon
\let\ex\exists
\let\fa\forall
\let\ps\supset
\let\dt\det

\def\rs#1{\raisebox{-0.4em}{$\big|_{#1}$}}

{\let\@noitemerr\relax
\vskip-2.7em\kern0pt\begin{abstract}
\noindent{\bf Abstract.}\es
We realize a given (monic) Alexander polynomial by a (fibered)
hyperbolic arborescent knot and link of any number of components,
and by infinitely many such links of at least 4 components. As
a consequence, a Mahler measure minimizing polynomial, if it
exists, is realized as the Alexander polynomial of a fibered
hyperbolic link of at least 2 components. For given polynomial,
we give also an upper bound on the minimal hyperbolic volume of
knots/links, and contrarily, construct knots of arbitrarily large
volume, which are arborescent, or have given free genus at least 2.
\\[2mm]
\end{abstract}
}


\section{Introduction}

The hyperbolic volume $\vol(L)$ of (the complement in $S^3$ of) a link
$L$ is an important, but not easy to understand geometric invariant.
For some time its relations to other topological and quantum
invariants, in particular the Alexander $\Dl$ \cite{Alexander}
and Jones \cite{Jones} polynomial, have been sought. In that
regard, recently a variety of connections between the hyperbolic
volume and "Jones-type" invariants has come to attention. Among
others, one would like to understand what geometric complexity
is measured by the polynomial invariants. An important question
remaining open is whether one can augment hyperbolic volume but 
preserve the Jones polynomial (or at least its Mahler measure).

In this paper, we will offer some analogous constructions for the
Alexander polynomial. First we realize the polynomial by a certain
arborescent knot. This yields an upper bound on the minimal volume
of a hyperbolic knot with given Alexander polynomial, which depends
only on the degree of the polynomial (see Theorem \ref{thbd}).
Apart from hyperbolicity, we show also that the knots have
canonical surfaces of minimal genus, and that these surfaces are
fiber surfaces if the Alexander polynomial is monic. Later we show how
to augment hyperbolic volume (Theorem \reference{theo2}). This first
construction simultaneously augments the slice genus. Another such
construction extends the result of Brittenham \cite{Brittenham2}.
It yields knots of arbitrarily large volume with given free
genus at least 2, with the additional feature that we can again
realize a given Alexander polynomial (Theorem \ref{thfg}).

A main theme will be to consider also various questions for
links. The realization result is extended first to links of
two (Theorem \reference{th2c}), and then of more components
(Theorem \reference{HHHH}). The hyperbolicity proof is,
unlike for knots, more involved, and requires the main effort.
It uses heavily the results of Oertel \cite{Oertel} and
Wu \cite{Wu}. A motivation was that for fibered links
of given polynomial not even primeness issues seem to have
ever been settled (and for more than 2 components, no
possibly prime links have been available). Another motivation,
and now application (Corollary \ref{cI}), is to confirm a
claim of Silver and Williams, that a polynomial of minimal
(positive) Mahler measure, if it exists, is realized as the
Alexander polynomial of a fibered hyperbolic 2-component link
(see Remark \reference{rI}).

Later we succeed in partially extending the construction to
obtain infinite families of links. An analogue of the infinite
realizability result of Morton \cite{Morton} for fibered knots
is shown for (arborescent) links of $\ge 4$ components (Proposition
\ref{i6,}), even for canonical fiber surfaces (for which it is
known not to hold in some other cases \cite{canon}).
Table \reference{tab2} at the end of the paper summarizes
these (in the context of some previous related) results.

We will use several methods, including Seifert matrices and skein
relations (for realizing Alexander polynomials), tangle surgeries
and Stallings twists (for generating infinite families of links),
some cut-and-paste arguments (for showing hyperbolicity),
and results of Gabai \cite{Gabai2,Gabai3} based on his
sutured manifold theory \cite{Gabai} (to prove fibering).

Constructions in a similar, but somewhat different, spirit were
proposed recently by Kalfagianni \cite{Kalfagianni}, Nakamura
\cite{Nakamura}, and Silver and Whitten \cite{SWh}. Most
properties studied there can be obtained from our work, too
(except for the knot group homomorphism in \cite{SWh}; see
remarks  \reference{rKf}, \reference{rSW} and \reference{rN}).
If one is mainly interested in Alexander polynomials and large
volume (but not in genera, fibering and arborescency), there are
generalizations in a further direction \cite{Friedl}, using
Kawauchi's imitation theory.

A worth remarking (though beyond our scope here) other connection 
is the Volume conjecture \cite{MM}, which asserts that one can 
determine (theoretically, but not practically) the volume \em{exactly}
from the Jones polynomial \em{and all} its cables. There is also
accumulating evidence that the (ordinary) Jones polynomial might
be able to provide in a different way (very practical) bounds on the
volume. Such bounds, which involve degrees or coefficients of the
polynomial, have been obtained for alternating links by Dasbach and 
Lin \cite{DL}, and later for Montesinos and 3-braid links by myself. 

On another related (but likewise not further pursued here) venue, I 
proved a conjecture of Dunfield \cite{Dunfield}, relating the 
determinant and
volume of alternating links. In particular, the determinant has an
exponential lower bound in terms of the volume. Since the determinant
can be expressed by both the Jones or Alexander polynomial,
we have a different relation of these invariants to the volume.
(Khovanov suggested a possible extension of Dunfield's correspondence
to non-alternating links, if instead of the determinant we take
the total dimension of his homology generalizing the Jones 
polynomial \cite{Khovanov}.)


\section{Some preliminaries\label{Sx}}

\subsection{Conway notation and Montesinos links\label{S21}}

\begin{defi}\label{dtg}
A \em{tangle} $Y$ is a set of two arcs and possible circles
(\em{closed components}) properly embedded in a ball $B(Y)$.
Tangles are considered up to homeomorphisms of $B(Y)$ that
keep fixed its boundary $\pa B(Y)$. Two tangles are \em{equivalent}
(in the sense of \cite{Wu}), if they are transformed by a
homeomorphism of their ball that preserves (but does not necessarily
fix) the 4 punctures of the boundary.
\end{defi}

\begin{figure}[htb]
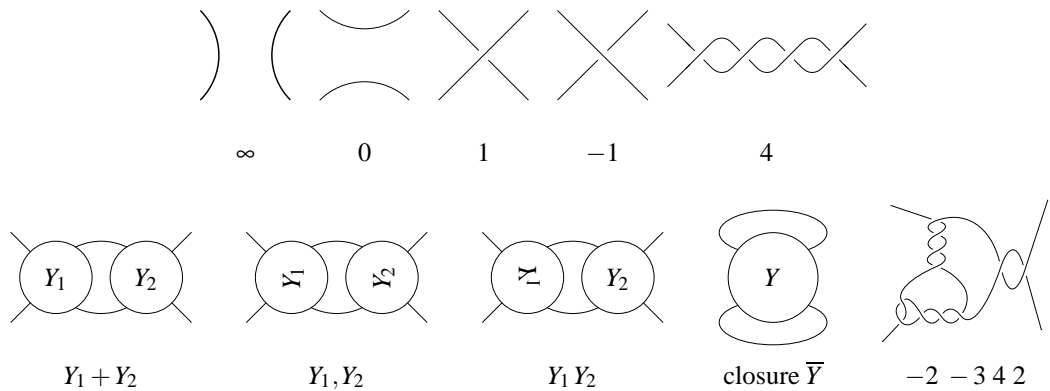

\[
\begin{array}{*5c}
\diag{6mm}{2}{2}{
  \piccirclearc{-1 1}{1.41}{-45 45}
  \piccirclearc{3 1}{1.41}{135 -135}
} & 
\diag{6mm}{2}{2}{
  \piclinewidth{50}
  \piccirclearc{1 -1}{1.41}{45 135}
  \piccirclearc{1 3}{1.41}{-135 -45}
} &
\diag{6mm}{2}{2}{
  \piclinewidth{50}
  \picmultiline{-15 1 -1.0 0}{2 0}{0 2}
  \picmultiline{-15 1 -1.0 0}{0 0}{2 2}
} &
\diag{6mm}{2}{2}{
  \piclinewidth{50}
  \picmultiline{-15 1 -1.0 0}{0 0}{2 2}
  \picmultiline{-15 1 -1.0 0}{2 0}{0 2}
} &
\diag{6mm}{4}{1}{
  \piclinewidth{50}
  \picmultigraphics{4}{1 0}{
    \picline{0.2 0.8}{0.8 0.2}
    \picmultiline{-14 1 -1.0 0}{0.2 0.2}{0.8 0.8}
  }
  \picmultigraphics{3}{1 0}{
    \piccirclearc{1 0.6}{0.28}{45 135}
    \piccirclearc{1 0.4}{0.28}{-135 -45}
  }
  \picline{-0.2 -0.2}{0.3 0.3}
  \picline{-0.2 1.2}{0.3 0.7}
  \picline{4.2 -0.2}{3.7 0.3}
  \picline{4.2 1.2}{3.7 0.7}
} \\[9mm]
\infty & 0 & 1 & -1 & 4 
\end{array}
\]
\[
\begin{array}{c*4{@{\qquad}c}}
\diag{6mm}{4}{2}{
  \piclinewidth{50}
  \picline{0 0}{0.5 0.5}
  \picline{0 2}{0.5 1.5}
  \picline{4 0}{3.5 0.5}
  \picline{4 2}{3.5 1.5}
  \piccirclearc{2 0.5}{1.3}{45 135}
  \piccirclearc{2 1.5}{1.3}{-135 -45}
  \pictranslate{1 1}{
    \picrotate{0}{
      \picscale{1 1}{
	\picfilledcircle{0 0}{0.8}{$Y_1$}
      }
    }
  }
  \pictranslate{3 1}{
    \picrotate{0}{
      \picscale{1 1}{
	\picfilledcircle{0 0}{0.8}{$Y_2$}
      }
    }
  }
} &
\diag{6mm}{4}{2}{
  \piclinewidth{50}
  \picline{0 0}{0.5 0.5}
  \picline{0 2}{0.5 1.5}
  \picline{4 0}{3.5 0.5}
  \picline{4 2}{3.5 1.5}
  \piccirclearc{2 0.5}{1.3}{45 135}
  \piccirclearc{2 1.5}{1.3}{-135 -45}
  \pictranslate{1 1}{
    \picrotate{90}{
      \picscale{1 1}{
	\picfilledcircle{0 0}{0.8}{$Y_1$}
      }
    }
  }
  \pictranslate{3 1}{
    \picrotate{90}{
      \picscale{1 1}{
	\picfilledcircle{0 0}{0.8}{$Y_2$}
      }
    }
  }
} &
\diag{6mm}{4}{2}{
  \piclinewidth{50}
  \picline{0 0}{0.5 0.5}
  \picline{0 2}{0.5 1.5}
  \picline{4 0}{3.5 0.5}
  \picline{4 2}{3.5 1.5}
  \piccirclearc{2 0.5}{1.3}{45 135}
  \piccirclearc{2 1.5}{1.3}{-135 -45}
  \pictranslate{1 1}{
    \picrotate{90}{
      \picscale{-1 1}{
	\picfilledcircle{0 0}{0.8}{$Y_1$}
      }
    }
  }
  \picfilledcircle{3 1}{0.8}{$Y_2$}
} &
\diag{6mm}{2.4}{3}{
  \piclinewidth{50}
  \picmultigraphics{2}{0 2}{
    \picellipse{1.2 0.5}{1.2 0.5}{}
  }
  \picfilledcircle{1.2 1.5}{1}{$Y$}
} & \eepsfs[2cm]{t1-rat234}
\\[9mm]
\text{ $Y_1+Y_2$ } & \text{ $Y_1,Y_2$ } & \text{ $Y_1\,Y_2$ } &
\text{ closure $\overline{Y}$ } & -2\ -3\ 4\ 2
\end{array}
\]
\caption{Conway's primitive tangles and operations with them. 
\label{fig1}}
\end{figure}

Figure \ref{fig1} shows the elementary tangles, tangle operations
and notation, mainly leaning on Conway \cite{Conway}. A \em{clasp}
is one of the elementary tangles $\pm 2$ and its rotations.
For two tangles $Y_1$ and $Y_2$ we write $Y_1+Y_2$
for the \em{tangle sum}. This is a tangle obtained by
identifying the NE end of $Y_1$ with the NW end of 
$Y_2$, and the SE end of $Y_1$ with the SW end of 
$Y_2$. The \em{closure} of a tangle $Y$ is a link
obtained by identifying the NE end of $Y$ with its 
NW end, and the SE end with the SW end. The
closure of $Y_1+Y_2$ is called \em{join} $Y_1\cup
Y_2$ of $Y_1$ and $Y_2$.

\begin{defi}
A link diagram is \em{arborescent}, if it can be obtained from the
tangles in figure \reference{fig1} by the operations shown therein.
An alternative description is as follows. Take a one crossing (unknot)
diagram. Repeat replacing some (single) crossing by a clasp (of any
orientation or sign). The diagrams obtained this way are exactly
the arborescent diagrams.
In Conway's \cite{Conway} terminology, these are diagrams with Conway
polyhedron $1^*$. A link is said to be arborescent if it admits an
arborescent diagram.

A graph $G$ is \em{series parallel}, if it can be obtained from 
\ $\diag{1em}{2}{1}{\llpoint{0 0.5}\llpoint{2 0.5}\picline{0 0.5}
{2 0.5}}$\ \ by repeated edge bisections and doublings. Such graphs 
correspond to arborescent link diagrams via the checkerboard graph
construction (see \cite{Kauffman,Mighton,Thistle} for example).
\end{defi}

\begin{defi}
A \em{rational} tangle diagram is the
one that can be obtained from the primitive Conway tangle
diagrams by iterated left-associative product in the way displayed
in figure \ref{fig1}. (A simple but typical example of
is shown in the figure.)
\end{defi}

\begin{figure}[htb]
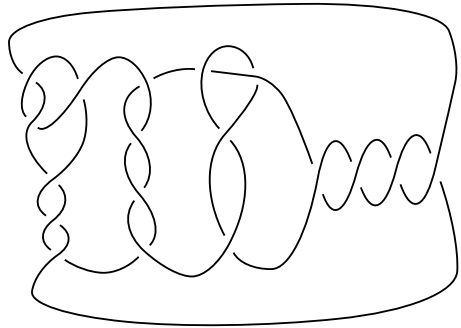

\[
\begin{array}{c}
\epsfsv{5cm}{t1-dioph3} \\
\end{array}
\]
\caption{The Montesinos knot \unhbox\@tempboxa\ with
Conway notation $(213,-4,22,40)$.\label{fig2}}
\end{figure}

Let the \em{continued} (or \em{iterated}) \em{fraction}
$[[s_1,\dots,s_r]]$ for integers $s_i$ be defined inductively by
$[[s]]=s$ and
\[
[[s_1,\dots,s_{r-1},s_r]]=s_r+\frac{1}{[[s_1,\dots,s_{r-1}]]}\,.
\]
The rational tangle $T(p/q)$ is the one with Conway notation
$c_{1}\ c_{2}\ \dots c_{n}$, when the $c_i$ are chosen so that
\begin{eqn}\label{ci}
[[c_{1},c_{2},c_{3},\dots,c_n]]=\frac{p}{q}\,.
\end{eqn}
One can assume without loss of generality that $(p,q)=1$, and $0<q<|p|$.
A \em{rational (or 2-bridge) link} $S(p,q)$ is the closure of $T(p/q)$.

Montesinos links (see e.g. \cite{BurZie}) are generalizations of
pretzel and rational links and special types of arborescent links. They
are denoted in the form $M(\frac{q_1}{p_1},\dots,\frac{q_n}{p_n},e)$,
where $e,p_i,q_i$ are integers, $(p_i,q_i)=1$ and $0<|q_i|<p_i$.
Sometimes $e$ is called the \em{integer part}, and the $\frac{q_i}
{p_i}$ are called \em{fractional parts}. They both together form the
\em{entries}. If $e=0$, it is omitted in the notation.

If all $|q_i|=1$, then the Montesinos link $M(\pm \frac{1}{p_1},\dots,
\pm \frac{1}{p_n},e)$ is called a \em{pretzel link}, of type $(\pm p_1,
\dots,\pm p_n,\eps,\dots,\eps)$, where $\eps=\sgn(e)$, and there
are $|e|$ copies of it.

To visualize the Montesinos link from a notation,
let $p_i/q_i$ be continued fractions of rational tangles
$c_{1,i}\dots c_{n_i,i}$ with $[[c_{1,i},c_{2,i},c_{3,i},\dots,
c_{l_i,i}]]=\frac{p_i}{q_i}$. Then $M(\frac{q_1}{p_1},
\dots,\frac{q_n}{p_n},e)$ is the link that corresponds to the
Conway notation
\begin{eqn}\label{cZ}
(c_{1,1}\dots c_{l_1,1}), (c_{1,2}\dots c_{l_2,2}), \dots,
(c_{1,n}\dots c_{l_n,n}), e\,0\,.
\end{eqn}
The defining convention is that all $q_i>0$ and if $p_i<0$, then
the tangle is composed so as to give a non-alternating sum with
a tangle with $p_{i\pm 1}>0$. This defines the diagram up to mirroring.
We sometimes denote the \em{Montesinos tangle} with Conway notation
\eqref{cZ} in the same way as its closure link.

An easy exercise shows that if $q_i>0$ resp. $q_i<0$, then 
\begin{eqn}\label{Mi}
M(\dots,q_i/p_i,\dots,e)\,=\,
M(\dots,(q_i\mp p_i)/p_i,\dots,e\pm 1)\,,
\end{eqn}
i.e. both forms represent the same link (up to mirroring).

Note that our notation \em{may differ} from other authors' by the sign
of $e$ and/or multiplicative inversion of the fractional parts. For
example $M(\frac{q_1}{p_1},\dots,\frac{q_n}{p_n},e)$ is denoted 
as ${\frak m}(e;\frac{p_1}{q_1},\dots,$ $\frac{p_n}{q_n})$ in
\cite[definition 12.28]{BurZie} and as $M(-e;(p_1,q_1),\dots,$ $
(p_n,q_n))$ and the tables of \cite{Kawauchi}.

Our convention chosen here appears more natural~-- the identity
\eqref{Mi} preserves the sum of all entries, and an integer entry
can be formally regarded as a fractional part. Theorem 12.29 in
\cite{BurZie} asserts that the entry sum, together with the vector
of the fractional parts, modulo $\bZ$ and up to cyclic permutations
and reversal, determine the isotopy class of a Montesinos link $L$.
So the number $n$ of fractional parts is an invariant of $L$;
we call it the \em{length} of $L$.

\setbox\@tempboxa=\hbox{$M(\myfrac{3}{11},-\myfrac{1}{4},
\myfrac{2}{5},4)$}

If the length $n<3$, an easy observation shows that the
Montesinos link is in fact a rational link. Then we could write
rational links as Montesinos links of length 1. For example,
$M(1)=M(\infty)$ is the unknot, and $M(0)$ is the 2-component
unlink, while $M(2/5)=M(5/2)$ is the figure-8 knot. This
simplification is not right, though, for Montesinos tangles
with $n=2$. Thus we keep (and will need) the length-2 notation
for tangles.


\subsection{Diagrams and geometric invariants}

\begin{defi}
A crossing in an oriented diagram looking like $\diag{7mm}{1}{1.2}{
  \pictranslate{0 0.1}{
  \picmultivecline{0.18 1 -1.0 0}{0 0}{1 1}
  \picmultivecline{0.18 1 -1.0 0}{0 1}{1 0}
  }
}$ is called \em{positive}, and $\diag{7mm}{1}{1.2}{
  \pictranslate{0 0.1}{
  \picmultivecline{0.18 1 -1.0 0}{0 1}{1 0}
  \picmultivecline{0.18 1 -1.0 0}{0 0}{1 1}
  }
}$ is a \em{negative} crossing. This dichotomy is called also \em{%
(skein) sign}. In an oriented diagram a clasp is called \em{positive},
\em{negative} or \em{trivial}, if both crossings are positive/negative,
resp. of different sign. Depending on the orientation of the involved
strands we distinguish between a \em{reverse
clasp} $ \diag{7mm}{2}{1.2}{
  \picPSgraphics{0 setlinecap}
  \pictranslate{1 0.6}{
    \picrotate{-90}{
      \rbraid{0 -0.5}{1 1}
      \rbraid{0 0.5}{1 1}
      \pictranslate{-0.5 0}{
      \picvecline{0.02 -.95}{0 -1}
      \picvecline{0.98 .95}{1 1}
    } } } }\es $
and a \em{parallel clasp}
$ \diag{7mm}{2}{1.2}{
  \picPSgraphics{0 setlinecap}
  \pictranslate{1 0.6}{
    \picrotate{-90}{
      \rbraid{0 -0.5}{1 1}
      \rbraid{0 0.5}{1 1}
      \pictranslate{-0.5 0}{
      \picvecline{0.02 .95}{0 1}
      \picvecline{0.98 .95}{1 1}
    } } } }\,.$
So a clasp is reverse if it contains a full Seifert circle,
and parallel otherwise. (We refer to \cite{Lickorish,Rolfsen} for
the notion of a Seifert circle.)
\end{defi}

For the later explanations, we must introduce the notion of twist
equivalence of crossings. The version of this relation we present
here follows its variants studied in \cite{gen1,gen2}. 

\begin{defi}\label{df22}
We say two crossings $p$ and $q$ of a diagram $D$ to be
\em{$\sim$-equivalent}, resp.\ \em{$\ssim$-equivalent}, if up
to flypes they form a reverse resp.\ parallel clasp. We remarked
in \cite{gen1} that $\sim$ and $\ssim$ are equivalence relations.
%
We write $t_{\sim}(D)$ for the number of or $\sim$-equivalence
classes of crossings in $D$. Set $t_{\sim}(K)$, the \em{reverse twist
number} of a knot or link $K$, to be the minimum of $t_{\sim}(D)$
taken over all diagrams $D$ of $K$.
%

In \cite{gen1} we noticed also that if
$p\sim q$ and $p\ssim r$, then $p=q$ or $p=r$. (There is the,
not further troubling however, exception that $D$ is the 2-crossing
Hopf link diagram, or has such a diagram occurring as a connected sum
factor.) So the relation $(p\sim q\lor p\ssim q)$ is also an
equivalence relation. We call this relation \em{twist equivalence}.
Thus two crossings are twist equivalent if up to flypes they form a
clasp. We will often call twist equivalence classes of crossings
in a diagram simply \em{twists}. (Some twists may consist of a single
crossing.) Let $t(D)$ denote the \em{twist number} of a diagram $D$,
which is the number of its twists. The twist number $t(K)$ of a knot
or link $K$ is the minimal twist number of any diagram $D$ of $K$.
Clearly $t(D)\le t_{\sim}(D)$ and $t(K)\le t_{\sim}(K)$. 
\end{defi}

With this terminology, we can state the following inequality we need:

\begin{theorem}(\cite{Lackenby}\label{LAT})
For a non-trivial diagram $D$ of a link $L$, we
have $10\,V_0\,\bigl(\,t(D)-1\,\bigr)\,\ge\,\vol(L)$\,,
where $V_0=\vol(4_1)/2\approx 1.01494$ is the volume of the ideal
tetrahedron.
\end{theorem}

Such an inequality, with the constant $10$ replaced by $16$, follows
from well-known facts about hyperbolic volume (see for example the
explanation of \cite{Brittenham}). Lackenby \cite{Lackenby} (whose
main merit is a lower volume bound for alternating links) repeated
this observation, and Agol-Thurston found, in the appendix to
Lackenby's paper, the optimal constant $10$, which is used below
for a better estimate.

\begin{rem}
Our notion of twist equivalence is slightly more relaxed than
what was called this way in \cite{Lackenby}, the difference
being that there flypes were not allowed. We call Lackenby's
equivalence here \em{strong twist equivalence}. However, it was
repeatedly observed that by flypes all twist equivalent crossings
can be made strongly twist equivalent, which Lackenby formulated
as the existence of \em{twist reduced} diagrams. Thus, assuming that
the diagram is twist reduced, we can work with twist equivalence
in our sense as with twist equivalence in Lackenby's sense (or
strong twist equivalence in our sense).
\end{rem}

A diagram is \em{special} if no Seifert circle contains
other Seifert circles in both interor and exterior.

\begin{defi}
A \em{Seifert surface} $S$ for an oriented link $L$ is a compact
oriented surface bounding $L$. A Seifert surface is \em{free}
if its complement is a handlebody. It is \em{canonical}, if it
is obtained by Seifert's algorithm from some diagram of $L$.
A \em{slice surface} is a surface properly embedded in $B^4$
whose boundary is $L\subset S^3$. We denote by
$g(L)$, $g_c(L)$, $g_f(L)$ and $g_s(L)$ the Seifert, canonical,
free and smooth slice genus of $L$. These are the minimal genera
of a (canonical/free) Seifert or slice surface of $L$, resp.
For a link $L$ we write $\chi(L)$, $\chi_c(L)$ and $\chi_s(L)$
for the analogous Euler characteristics (we will not need $\chi_f$).
\end{defi}

Seifert's algorithm is explained, for example, in \cite{Rolfsen}.
We will use also some of the detailed discussion given to it in
\cite{gen1,gen2}.

A canonical Seifert surface is free, and any Seifert surface is a
slice surface. Thus $g_s(K)\le g(K)\le g_f(K)\le g_c(K)$ for any
knot $K$. By $u(K)$ we denote the \em{unknotting number} of $K$. 
Then it is known that $g_s(K)\le u(K)$.

For a link $L$, let $n(L)$ be the \em{number of
components} of $L$. Then $\chi_{[s/c]}(L)=2-n(L)-2g_{[s/c]}(L)$.

\subsection{The Alexander-Conway polynomial}

\begin{defi}\label{dad}
Below it will be often convenient to work with the \em{Conway
polynomial} $\nb(z)$. It is given by the value $1$ on the
unknot and the \em{skein relation}
\begin{eqn}\label{sRs}
\nb(D_+)-\nb(D_-)=z\nb(D_0)\,.
\end{eqn}
Here $D_{\pm}$ are diagrams differing only at one crossing,
which is positive/negative, and $D_0$ is obtained by smoothing
out this crossing. The Conway polynomial is equivalent to
the (1-variable\footnote{In this paper Alexander polynomials
are always understood to be the 1-variable versions.}) Alexander
polynomial $\Dl$ by the change of variable:
\begin{eqn}\label{nbDl}
\nb(t^{1/2}-t^{-1/2})=\Dl(t)\,.
\end{eqn}
For that reason we will feel free to exchange one polynomial
for the other whenever we deem it convenient. For knots
$\nb\in 1+z^2\bZ[z^2]$ and for $n$-component
links (with $n>1$) we have $\nb\in z^{n-1}\bZ[z^2]$.
We call such $\nb$ and the corresponding $\Dl$ \em{admissible}
polynomials. Each admissible polynomial is indeed realized
by some knot or link.
\end{defi}

There is another description for $\Dl$. Given a Seifert surface $S$
of genus $n=g(S)$ for a knot $K$, one associates to it a \em{Seifert
matrix} $V$ (a $2n\times 2n$ matrix of integer coefficients), and
we have
\[
\Dl(t)=t^{-n}\det (V-tV^T)\,,
\]
where $V^T$ is the transposed of $V$. This is described in
\cite{Rolfsen}, for example. 

A direct understanding of the relation between the skein-theoretic
and Seifert-matrix-related properties of $\Dl$ is still a
major mystery in knot theory. Solving it may shed light on a
topological meaning of the newer polynomials. To the contrary,
the long-term lack of such a meaning justifies the pessimism
in expecting the desired relation. Nonetheless, both descriptions
of $\Dl$ offer two independent ways of keeping control on it, and
we will successfully combine them in some of the below constructions.

We remark also that $\nb$ (and $\Dl$) is symmetric resp.
antisymmetric w.r.t. taking the mirror image, depending on the
odd resp. even parity of the number of components. This means
in particular that amphicheiral links of an even number of
components have vanishing polynomial. (Here amphicheirality
means that an isotopy to the mirror image is to preserve or reverse
the orientation of \em{all} components \em{simultaneously}, while
it is allowed components to be permuted.)

\begin{defi}\label{defmon}
Let $[X]_{t^a}=[X]_a$ be the coefficient of $t^a$ in a polynomial
$X\in\bZ[t^{\pm1}]$.  For $X\ne 0$, let 
$\cC_X\,=\,\{\,a\in\bZ\,:\,[X]_a\ne 0\,\}$ and 
\[
\md X=\min\,\cC_X\,,\quad \Md X=\max\,\cC_X\,,
\mbox{\qquad and\qquad} \spn X=\Md X-\md X\, 
\]   
be the \em{minimal} and \em{maximal degree} and \em{span} (or 
breadth) of $X$, respectively. The \em{leading coefficient}
$[X]_{*}$ of $X$ is defined to be $[X]_{\Md X}$. If this
coefficient is $\pm 1$, we call $X$ \em{monic}.
\end{defi}

A link in $S^3$ is \em{fibered} if its complement is a surface bundle
over $S^1$. By a classical theorem of Neuwirth-Stallings, the
fiber is then a minimal genus Seifert surface, and such a Seifert
surface is unique. The operations \em{Hopf (de)plumbing} and
\em{Stallings twist} are described, for example, in Harer
\cite{Harer}. (A Stallings twist is a $\pm 1$ surgery along an
unknot in the complement of the fiber surface, which can be
isotoped into the fiber.) Harer showed
that every fiber surface in $S^3$ can be constructed
from a disk by a sequence of these operations. Besides, there is
Gabai's geometric work to detect (non-)fiberedness \cite{Gabai4}.
We call a fibered link $L$ \em{canonically fibered} if its fiber
surface can be obtained by Seifert's algorithm on some diagram of $L$.

It is known that $\Md\Dl(K)\le g(K)$ for any knot $K$, and
similarly $2\Md\Dl(L)\le 1-\chi(L)$ for any link $L$.
The Alexander polynomial of a \em{fibered} link $L$ satisfies
$2\Md\Dl(L)=1-\chi(L)$ and $[\Dl]_*=\pm 1$ (see \cite{Rolfsen}).

By $\br{x}$ we will mean the greatest integer not greater than $x$,
and $\BR{x}$ denotes the smallest integer not smaller than $x$.

\section{Small volume knots}

In this section we will consider the problem how one can estimate
the volume of a hyperbolic knot in terms of the Alexander polynomial.
Simultaneously, we will try to estimate
the various genera (and for links, Euler characteristics).
For instance, it makes sense to ask 

\begin{question}\label{mvt}
What is the minimal twist number, or the minimal volume of a
hyperbolic knot, with given Alexander polynomial?
\end{question}

As the Alexander polynomial provides upper bounds on the crossing
number of \em{alternating} knots \cite{Crowell}, it certainly does
so for the twist number (and volume). Dunfield's correspondence
mentioned in the introduction is a sharper version of this easy
observation.
There exist also, for arbitrary knots, lower bounds on the twist
number from the Alexander polynomial, as we prove in
joint work with Dan Silver and Susan Williams \cite{STW}.

Note that one must exclude non-hyperbolic knots if we consider
the volume in question \ref{mvt}. Otherwise take a knot $K$ 
realizing $\Dl$. Then a satellite around $K$ with an 
unknotted pattern of algebraic degree 1, but geometric
degree $>1$, has the same Alexander polynomial.

%

The following result gives some information on
question \ref{mvt}.

\begin{theorem}\label{thbd}
Assume $\Dl\in\bZ[t^{\pm 1}]$ satisfies let $\Dl(t)=\Dl(1/t)$,
$\Dl(1)=1$, and let $\Md\Dl=d$. Then there is an arborescent knot
$K$ with the following properties.
\begin{enumerate}
\item \label{prtA} We have $\Dl(K)=\Dl$, $u(K)\le 1$, and
  $t_{\sim}(K)\le 4d-1$ if $d>0$.
\item \label{prtA'} A Seifert surface $S$ of genus $d$ for $K$ is
  obtained as a canonical surface of a special arborescent diagram
  of $K$. In particular $g(K)=g_c(K)=d$, so $S$ is of minimal genus.
\item \label{prtC} If $\Dl$ is monic, then $S$ is a fiber surface.
\item \label{prtB} If $\Dl$ is not the unknot or trefoil
  polynomial, then $K$ is 
  hyperbolic, and
  \begin{eqn}\label{vbd}
    0<\vol(K)\le 10V_0(4d-3)\,.
  \end{eqn}
\end{enumerate}
\end{theorem}

\begin{rem}
By a result of Hirasawa \cite{Hirasawa}, a canonical surface
from some diagram $D$ of a link $L$ is always canonical
w.r.t. a special diagram $D'$ of $L$. However, the procedure
he uses to turn $D$ into $D'$ does not preserve
arborescency (of the diagram).
\end{rem}

\begin{rem}\label{zas}
It follows from \cite{Kobayashi,gen1} that another knot of $g_c,u\le
1$ cannot have the Alexander polynomial of the unknot or trefoil.
Contrarily, if we waive on $u\le 1$ (and on fibering, and $g_c(K)=0$
for $\Dl=1$), then there is an infinity of pretzel knots
$(p,q,r)$ for $p,q,r$ odd with such polynomials.
\end{rem}

\begin{exam}
Among trivial polynomial knots, the two 11 crossing knots are 
arborescent, of unknotting number one, and have $\vol\approx 
11.2$. The smallest volume knot with trivial polynomial
I found is the $(-3,5,7)$-pretzel knot, where $\vol\approx
8.5$, but it is not of unknotting number one.

The knot $13_{5111}$ of \cite{KnotScape} is arborescent, has
$u(K)=1$ and the trefoil polynomial, and $\vol\approx 11.3$.
The smallest volume knot I found with this polynomial is
$13_{8541}$ with $\vol\approx 7.8$, but it is (apparently)
not arborescent nor of unknotting number one.
\end{exam}

There have been several other previous constructions of (fibered) 
knots (and links) with given (monic) polynomial, for example 
\cite{Burde,Kanenobu,Levine,Morton,Quach}. The new main features
here are the volume estimate and arborescency and to somewhat
smaller extent genus minimality of the canonical surface.

\begin{rem}\label{rN}
A triggering point for the present work was Nakamura's study of
braidzel surfaces \cite{Nakamura3}. Using these, he showed in
\cite{Nakamura} that one can choose $K$ in part \reference{prtA}
of Theorem \reference{thbd}, so that it has braidzel genus $n$
(and unknotting number one), by realizing a Seifert matrix in
\cite{Seifert}. But these braidzel surfaces are unlikely canonical.
Then, simultaneously to this writing, he used a Seifert matrix
of Tsutsumi and Yamada \cite{TY} (see the below proof), to find
braidzel surfaces isotopic to canonical surfaces of $4d-1$ twists
\cite{Nakamura2}. (I was pointed to this matrix also by him;
previously I used the one he gave in \cite{Nakamura} with a
weaker outcome.) Thus he gives a method that combines all our
properties except hyperbolicity and arborescency. 
\end{rem}

A different construction, producing (arguably always)
hyperbolic knots, is due to Fujii \cite{Fujii}. His
knots have tunnel number one, and are 3-bridge, but are unlikely
arborescent, and do not (at least in an obvious way) realize the
canonical genus by the degree of $\Dl$. His diagrams have
unbounded twist number even for fixed degree, and 
a similar volume bound using Thurston's surgery theorem
appears possible, but more elaborate and likely less
economical than ours.

After finishing this work, we learned that the same knots were
considered by H.~Murakami in \cite{Murakami}. We will nonetheless
go beyond the reproduction of his result (which he uses with a
different motivation from ours) that these knots have the
proper Alexander polynomial.

\proof[of Theorem \reference{thbd}]
{\em parts \reference{prtA} and \reference{prtA'}.} 
Let $\nb(z)$ be the Conway version of $\Dl$, and
\begin{eqn}\label{16}
\nb(z)=\,1-a_1z^2+a_2z^4-a_3z^6+\dots\,+(-1)^{n}
a_dz^{2d}\,\in\,\bZ[z^2]\,,
\end{eqn}
for integers $a_1,\dots,a_d$, so
$a_i=(-1)^i[\nb]_{2i}$. 
By Tsutsumi and Yamada 
\cite{TY}, it suffices to realize the matrices $V_n$ (shown for $d=2,4$,
with omitted entries understood to be zero, and with the obvious
generalization to arbitrary $d$)
\begin{eqn}\label{57}
V_2\,=\,\left(
\begin{array}{cc|cc}
-1   & -1   &     &     \\
0    & a_1  & 1   &     \\
\hline
     & 1    & 0   & -1  \\
     &      & 0   & a_2 \\
\end{array}
\,\right)\,,
\quad
V_4\,=\,\left(
\begin{array}{cc|cc|cc|cc}
-1   & -1   &     &     &     &     &     &      \\
0    & a_1  & 1   &     &     &     &     &      \\
\hline
     & 1    & 0   & -1  &     &     &     &      \\
     &      & 0   & a_2 & 1   &     &     &      \\
\hline
     &      &     & 1   & 0   & -1  &     &      \\
     &      &     &     & 0   & a_3 & 1   &      \\
\hline
     &      &     &     &     & 1   & 0   & -1   \\
     &      &     &     &     &     & 0   & a_4  
\end{array}
\,\right)
\end{eqn}
as Seifert matrices of canonical surfaces. Then
$\Dl(t)=t^{-n}\det (V_n-tV_n^T)$.

The solution is given by a sequence of graphs. We display
the first three in figure \reference{figgr}; the series is
continuable in the obvious way. The example for $d=a_1=3$, 
$a_2=a_3=-2$ is shown as a knot diagram on the left side below. 

\begin{figure}[htb]
\[
\diag{1cm}{2}{2}{
  \lpoint{0 1}
  \lpoint{2 1}
  \piccurve{0 1}{0.6 0.0}{1.4 0.0}{2 1}
  \piccurve{0 1}{0.6 2.0}{1.4 2.0}{2 1}
  \picline{0 1}{2 1}
    \picputtext{1 1.93}{$-1$}
    \picputtext[u]{1 0.91}{$-1$}
  \picputtext{0.5 0.07}{$2a_1+1$}
}\rx{1.3cm}
{
\diag{2cm}{2.5}{1.7}{
  \piclinewidth{60}
  \pictranslate{0 0.1}{
    \point{0 1}
    \point{2 1}
    \point{1.7 0.72}
    \piccurve{0 1}{0.6 0.3}{1.4 0.3}{2 1}
    \piccurve{0 1}{0.6 1.7}{1.4 1.7}{2 1}
    \picline{0 1}{2 1}
    \piccurve{0 1}{0.2 -0.0}{1.4 -.2}{1.7 0.72}
    \picputtext{1 1.67}{$-1$}
    \picputtext[u]{1 0.95}{$-1$}
    \picputtext[l]{1.9 0.8}{$2a_1-1$}
    \picputtext[l]{0.8 0.6}{$2$}
    \picputtext[l]{0.8 0.25}{$-1$}
    \picputtext[l]{1.65 0.45}{$-1$}
    \point{1.5 0.35}
    \piccurve{0 1}{0.2 -0.4}{1.2 -0.4}{1.5 0.35}
    \picputtext[l]{1.5 0.1}{$2a_2+1$}
  }
}
}\rx{1.2cm}
\diag{2.2cm}{2.5}{2.3}{
  \piclinewidth{60}
  \pictranslate{0 0.6}{
    \point{0 1}
    \point{2 1}
    \point{1.7 0.72}
    \piccurve{0 1}{0.6 0.3}{1.4 0.3}{2 1}
    \piccurve{0 1}{0.6 1.7}{1.4 1.7}{2 1}
    \picline{0 1}{2 1}
    \piccurve{0 1}{0.2 -0.0}{1.4 -.2}{1.7 0.72}
    \picputtext{1 1.67}{$-1$}
    \picputtext[u]{1 0.95}{$-1$}
    \picputtext[l]{1.9 0.8}{$2a_1-1$}
    \picputtext[l]{0.8 0.6}{$2$}
    \picputtext[l]{0.8 0.25}{$-1$}
    \picputtext[l]{1.65 0.45}{$-1$}
    \point{1.5 0.35}
    \piccurve{0 1}{0.2 -0.4}{1.2 -0.4}{1.5 0.35}
    \picputtext[l]{1.5 0.1}{$2a_2-1$}
    \point{1.3 0.05}
    \piccurve{0 1}{0.1 -0.6}{1.1 -0.7}{1.3 0.05}
    \point{1.1 -0.3}
    \piccurve{0 1}{-.0 -0.8}{0.9 -0.9}{1.1 -0.3}
    \picputtext[l]{1.3 -0.15}{$-1$}
    \picputtext[l]{1.1 -0.45}{$2a_3+1$}
    \picputtext[l]{0.7 0.0}{$2$}
    \picputtext[l]{0.7 -.3}{$-1$}
    {
    \picputtext[l]{0.8 1.25}{$A$}
    \picputtext[l]{0.55 0.8}{$B$}
    \picputtext[l]{0.5 0.4}{$C$}
    \picputtext[l]{0.43  0.13}{$D$}
    \picputtext[l]{0.42 -0.12}{$E$}
    \picputtext[l]{0.36 -0.35}{$F$}
    }
  }
}
\]
\vspace{3mm}
\[
\vcbox{\epsfdiag{\epsfs{5cm}{t1-alrealiz}}{2.2cm}{
\picputtext{0.18 0.47}{*}
\picputtext{0.54 1.35}{**}
}}
\kern2cm
\vcbox{\epsfdiag{\epsfs{5cm}{t1-alrealiz2}}{1.1cm}{
  \picputtext{0.7 2.5}{$D$}
  \picputtext{2.6 2.4}{$C$}
  \picputtext{2.2 1.2}{$B$}
  \picputtext{4.2 4.3}{$A$}
}}
\]
\caption{\label{figgr}}
\end{figure}

One obtains the surfaces from the graphs as follows.  Each vertex
corresponds to a Seifert circle of valence $\ge 3$.
(The valence of a Seifert circle is the number of crossings
attached to it.)  Each edge with
label $x$ corresponds to a band of $|x|$ reverse half-turns of (skein)
sign $\sgn(x)$, enclosing $|x|-1$ valence-2-Seifert circles in between.

To obtain the Seifert matrix, for each of the bounded regions of the
complement of the graph, choose a loop going around the boundary.
The rows of $V_n$ (from top to bottom) and columns (from left to
right) correspond to loops ordered alphabetically by the letter
in their region. The orientation is coherently chosen, so two
loops pass along a common edge (twisted band) in opposite direction.
If the label of an inner edge is odd (always $-1$), the loops
are intertwined. Let them intersect once on one of the neighbored
Seifert circles, so as to reinstall their position.  Otherwise
loops do not intersect.

The graphs are series parallel (as defined in \S\ref{Sx}) so the
knots are arborescent.  Unknotting number one is visualized by drawing
the knot diagram. Resolving the parallel clasp * (the double edge
labeled $-1$ in the graph) gives an unknotting crossing change. 

{\em part \reference{prtC}.} 
Assume $\Dl$ is monic. We show that $S$ can be constructed from a
genus one fiber surface $S'$ by Hopf plumbings and Stallings twists.
To that vein, we apply them in reverse order and reduce $S$ to $S'$.

Deplumbing a Hopf band, one resolves one of the crossings in the
clasp * in the diagram on figure \reference{figgr}. A Hopf
(de)plumbing preserves the fiber property by \cite{Gabai2,Gabai3}.
By a Stallings twist, one cancels the other crossing, together 
with the twist of $2a_1-1$. Then one removes the Hopf link as 
connected sum factor (the clasp **) by deplumbing another Hopf
band. By iterating this procedure, one reduces $K$ to a diagram of
a negative clasp and a twist of $+3$ or $-1$ (since $\Dl$ is monic).
This is the fiber surface $S'$ of the trefoil or figure-eight knot.

{\em part \reference{prtB}.} 
By the work of Hatcher and Thurston, we must argue that the
knots are not satellite, composite or torus knots. It is known
from \cite{Oertel,Wu2} that arborescent knot complements are
atoroidal, so there is no satellite or composite arborescent
knot. Arborescent torus knots are classified in the monograph
of Bonahon-Siebenmann \cite{BS}, which is only told to exist.
However, we can use a published argument. In our case also
$u(K)=1$, and only the trefoil is a torus knot of unknotting
number one. This probably first follows from the signature formulas
of torus knots \cite{GLM,Hirzebruch}, or more directly from the
subsequent result of \cite{KM}. So we have hyperbolic knots $K$ except
the trefoil and unknot.

We have a diagram $D$ with $4d-1$ $\sim$-equivalence classes,
with two of them (of a single crossing each; the boundary of region
$A$) forming the parallel clasp *, so $4d-2$ twist equivalence classes. 
Then applying Theorem \reference{LAT}, we have the
stated volume estimate. \qed
%

\begin{rem}
For an infinite series of knots, we can apply tangle surgery (see
below), at the cost of slightly increasing the twist number. (However, 
it is not evident how to preserve fiberedness; see the remarks in
\S\ref{S72}.)
\end{rem}

\section{Two component links\label{S4}}

With some more work, we can obtain a result of
almost the same stature as Theorem \ref{thbd}
for links of two components.

\begin{theorem}\label{th2c}
Any admissible Alexander polynomial of a 2-component link
is realized by an arborescent link $L$ with $d=2\Md\Dl= 1-\chi_c(L)$,
which can be chosen to have the following further properties.
\begin{enumerate}
\item If $\Dl$ is monic, then $L$ is additionally fibered.
\item If $d>1$ (that is, $\nb(z)\ne kz$, $k\in \bZ$), then $L$
is hyperbolic, and
\begin{eqn}\label{vbd2}
0<\vol(L)\le 20V_0(d-1)\,.
\end{eqn}
\end{enumerate}
\end{theorem}

\begin{rem}\label{rI}
Silver and Williams were interested in proving, that if Lehmer's
question on the existence of a Mahler measure minimizing 
polynomial $f$ has an affirmative answer, then $f$ can be chosen
to be the Alexander polynomial of a fibered hyperbolic knot
or 2-component link. They claimed this in a preliminary (arXiv
v1) version of \cite{SW}, but there was an error in their
reasoning (as has been noted in the revision). The provision of a
correction motivated the study of two component links here.
However, this correction requires some work, as a ``pre-prepared''
argument, like in the case of knots, does not seem available.
\end{rem}

Theorem \reference{th2c}, beside confirming their claim, shows
a bit more. While it is of course more interesting if one can
exclude the 2-component links (or relatedly, to understand the
significance of the condition $\Dl(1)=1$ in Lehmer's question),
once links come in, our theorem first eliminates the (need of)
knots. We will see later, with Theorem \reference{HHHH}, that 
we can choose the number of link components arbitrarily (as long
as above $1$).

\begin{corr}\label{cI}
A polynomial of minimal Mahler measure (if such exists) is realized
as the Alexander polynomial of a fibered hyperbolic arborescent
2-component link. \qed
\end{corr}

However, second, we see that, from the point of view of 
mere realizability, there is nothing special to Lehmer's (or any 
other monic reciprocal) polynomial. This should caution in seeking 
a topological meaning behind Lehmer's question along these lines.

\proof[of Theorem \reference{th2c}]
To obtain a link $L$ of two components with given $\nb$,
smooth out the unknotting crossing in the knot found for
$1+z\nb$ in the proof of Theorem \ref{thbd}. Observe that
on the surface this is a Hopf deplumbing, so that fiberedness
is preserved for monic polynomials. The Conway polynomial is
$a_1z-a_2z^3+a_3z^5-\dots$, with the $a_i$ as in \eqref{16}.

The inequality \eqref{vbd2} is clear once we show hyperbolicity.
For this we assume that $a_1\not\in\{1,2,3\}$. Otherwise,
realize $-\Dl$, and take the mirror image.

We show below in Lemma \ref{lato} that $L$ is atoroidal
if $a_1\ne 1$. Atoroidality settled, hyperbolicity follows
from Hatcher-Thurston once Seifert fibred link complements
are excluded. Links with Seifert fibred complements are
determined by Burde and Murasugi
\cite{BM}. It follows from their work that all components
of such links are (possibly unknotted) torus knots.
Excluding the case of $d=1$, giving the $(2,\,.\,)$-torus
links, in our examples we have an (obviously) unknotted component
$O$, and a further component $K$. Now note that the knot $K$
is of the form that is obtained by our previous construction in
Theorem \reference{thbd}. By that construction,
\begin{eqn}\label{bnK1}
\nb_K\ne 1\,,
\end{eqn}
so $K$ is knotted. Also, by the proof of part \ref{prtB}
of that theorem, $K$ is hyperbolic (and in particular
not a torus knot), unless it is a trefoil.
If $K$ is a trefoil, the proof in \cite{BM}
shows that a 2-component link of an unknot and a trefoil
occurs only in their case (b). A look at the argument there
shows that we must have $a_1=lk(K,O)\in\{\pm 2,\pm 3\}$.
This leaves only 4 links; they can be specified (up
to component orientation) as the closures of the 3-braids 
$\sg_1^{-2}\sg_2^{-2}\sg_1^{-1}\sg_2^{2-2a_1}$. A check with 
Jeff Weeks' software SnapPea, available as a part of 
\cite{KnotScape}, shows that for $a_1=-2,-3$ the links
are hyperbolic (while for $a_1=2,3$ they are not, which
explains the other initial restriction). \qed

\begin{defi}\label{dtw}
In the following a \em{twist of $x$} for $x\in \bZ$ is
understood to mean a twist of $|x|$ crossings of (skein) sign
$\sgn(x)$. We call $|x|$ the \em{length} of the twist. A twist
is \em{reverse} or \em{parallel} if the crossings it contains
are $\sim$ or $\ssim$-equivalent resp., according to definition
\reference{df22}. (A twist of a single crossing is
simultaneously both reverse and parallel.)
\end{defi}

In order to avoid that the 2-component link is a connected sum
with a Hopf link factor, we need $a_1\ne 1$. First, we prove

\begin{lemma} \label{prl}
The above constructed link $L$ is prime if $a_1\ne 1$.
\end{lemma} 

\proof An easy ``proof'' is a routine application of the technique
in \cite{KL}, but here is another proof (with a fully different
argument, and worth dropping the quotes).

Since we assume $d>1$, our link $L$ consists of an
(obviously) unknotted component $O$, and another component
$K$. We observed that $K$ actually is of the form that
was constructed in Theorem \reference{thbd}. Then we have
\eqref{bnK1}, so in particular $K$ is knotted. Moreover
\begin{eqn}\label{gf}
\Md\nb_K=\Md\nb_L-1\quad\mbox{and}\quad[\nb_K]_*= \pm [\nb_L]_*\,.
\end{eqn}
We also have $u(K)=1$, so that $K$ is prime by \cite{Scharl}.
Hence the only possible way that $L$ is composite is that $L=K\#
L'$, where $L'$ is a link of two unknotted components. Because of
\eqref{gf} we have $\nb_{L'}=\pm z$. By additivity of the genus
under connected sum, $L'$ must bound an annulus, and then, since
its both components are unknotted, $L'$ must be a Hopf link. Now
\[
a_1=lk(K,O)=[\nb_L]_z=\pm [\nb_K]_{z^0}=\pm 1\,.
\]
Since we excluded $a_1=1$, the sign is negative, and so $L'$ is a
negative Hopf link. Let $\wt L$ be the link obtained from $L$ 
by reversing the orientation of $O$.
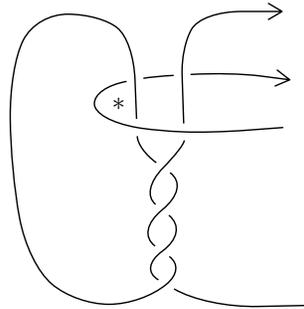
\begin{eqn}\label{gj}
\diag{1cm}{4}{4}{
  \picputtext[dl]{0 0}{\epsfs{4cm}{t1-KOlink}}
  \picputtext{1.45 2.71}{$*$}
}
\end{eqn}
Then we must have 
\begin{eqn}\label{st*}
\nb_{\wt L}=-\nb_L=+z\nb_K\,.
\end{eqn}
To show that this is not the case, we calculate $\nb_{\wt L}$.
Apply the skein relation \eqref{sRs} at the clasp $*$. (In $\wt L$
the orientation is so that the clasp is negative and parallel.)
\[
\nb(D_-)\,=\,\nb(D_+)\,-\,z\nb(D_0)\,.
\]
Then $D_+$ depicts the connected sum of a parallel $(2,4)$-torus
link with $K$, so $\nb(D_+)=(2z+z^3)\nb_K$. The diagram $D_0$
depicts a knot $K'$, which is obtained from $K$ by reversing the
sign of the crossings in the unknotting (parallel) clasp.

If $\nb_i$ are the polynomials of links $L_i$ with diagrams
equal except at one spot, where a parallel twist of $i$ positive
crossings is inserted, then by the skein relation
\[
\nb_4\,=\,\nb_2+z\nb_3\,=\,\nb_2+z\nb_1+z^2\nb_2\,=\,\nb_2+\nb_2-\nb_0
+z^2\nb_2\,=\,(2+z^2)\nb_2-\nb_0\,.
\]
So $\nb(D_0)=\nb(K')=z^2+2-\nb_K$. Then using \eqref{bnK1}, we have
\[
\nb_{\wt L}=\nb_{D_-}=(2z+z^3)\nb_K-2z-z^3+z\nb_K\,\ne\,z\nb_K\,,
\]
with the desired contradiction to \eqref{st*}. \qed

\begin{lemma}\label{lato}
The link $L$ is atoroidal if $a_1\ne 1$.
\end{lemma}

For the proof we require some cut-and-paste arguments. We lean
closely on the work of Wu \cite{Wu}. Let us fix some notation
and terminology first. All manifolds are assumed in general position,
so intersections are transversal. We use the formalism of tangle
operations in figure \ref{fig1} (see also the related explanation in 
and after Definition \reference{dtg}).

Writing again by $B(Y)$ the ball in which a tangle $Y$ lives,
we denote by $B(Y)\sm Y=X(Y)$ the \em{tangle space} of $Y$.
(This is a 3-manifold with a genus two surface as boundary;
see \cite{Wu}.) By $E(L)=S^3\sm L$ we denote the complement
of the link $L$.

We call a disk properly embedded in $X(Y)$ \em{separating}
if both balls in its complement contain parts of $Y$. We call a
tangle $Y$ \em{prime} \cite{KL} if it has no separating disk and
every sphere in $B(Y)$ intersecting $Y$ in two points bounds a
ball in $B(Y)$ intersecting $Y$ in an unknotted arc.

\proof[of Lemma \reference{lato}]
If $d=3$, then we have the Montesinos link
$M(-\frac{2a_2}{4a_2+1}, \frac{1}{2},\frac{1}{2a_1-2})$. 
Atoroidality follows then from \cite{Oertel}.
Our form is not among those given in corollary 5 there
(see in particular the proof of the corollary\footnote{but
beware that the hyperbolicity argument~-- which we do
not require~-- contains an error; see the remarks at the end of
\S\ref{S51} below.}). 

Let now $d\ge 5$. In our situation, $L=Y_1 \cup Y_2$ is a 2-component
link, and for integers $k\ne 0$, and $m$ odd we can write in the
notation of figure \reference{fig1}
\begin{eqn}\label{Tds}
Y_1\ =\ (U\ 1\ 1,-2)\ m\,,\quad\mbox{and}\quad Y_2\ =\ (2k,-2)\ 1\ 1\ \ 
(\mbox{$= R[2k,-2;1] + 1$ in the notation of \cite{Wu}})\,.
\end{eqn}
($U$ is a, possibly rational, arborescent tangle; $Y_2$ is the tangle
in \eqref{gj}.) So $Y_2$ has an unknotted closed component $O$, but
$Y_1$ has none. Let $K$ be the other, knotted, component of $L$.
It is easily verified using \cite{KL} that $Y_i$ are prime.

So now assume $T\subset E(L)$ is an essential (i.e., incompressible
and not boundary parallel) torus. $T$ bounds a solid torus $S$ we
call also \em{interior} $\ir T$. If $T$ bounds two solid tori,
$T$ is unknotted. Then choose one solid torus to be $S$. Let
$R=S^3\sm S$ be the other complementary \em{region}, which we call
also \em{exterior} $\et T$. Let $B_i=B(Y_i)$ be the balls in
which $Y_i$ are contained (with $B_1\cup B_2=S^3$), $X_i=X(Y_i)$
be the tangle spaces and $P=\bd X_1\cap \bd X_2$ their
common boundary, a 4-punctured sphere $C=\bd B_i$. 
We call $T$ \em{separating} if both regions of
$S^3\sm T$ contain one component of $L$ each.

\begin{slemma}\label{sl40}
Let $F\subset T$ be a circle, and assume $F$ bounds a disk $D$ in
one of the complementary regions of $T$, and $D$ is not parallel to
$T$. Then $|D\cap L|\ge 2$.
\end{slemma}

\proof An empty intersection is clearly out because $T$ is
incompressible. Assume $|D\cap L|=1$. We produce a contradiction
in cases by assuming that some meridional disk $D$ of $T$ intersects
$L$ in one point. (We choose the interior $S$ of $T$ to contain $D$.)

\begin{caselist}
\case $T$ is knotted.

\begin{caselist}
\case If $T$ is separating, the component $M$ of $L$ in $S=\ir
T$ is composite (and $T$ is a swallow torus) or satellite, or
$T$ is $\bd$-parallel to $M$. Now neither of the components
of $L$ is a composite or satellite knot (see proof of Theorem
\reference{thbd}, part \reference{prtB}), and $T$ is essential,
so we have a contradiction to all options.

\case If $T$ is not separating, then $L$ is the connected sum of
the knot type of $T$ with some 2-component link (obtained by
reembedding unknottedly $S=\ir T$). This contradicts
Lemma \reference{prl}.
\end{caselist}

\case So now consider the case $T$ is unknotted. Then $T$ must be
separating (otherwise it compresses in its exterior). But then
if $T$ is not $\bd$-parallel, then $L$ is the connected sum
of the component of $L$ in $S$ with a satellite of the Hopf link
(with a pattern that keeps the core of $S$). This again
contradicts Lemma \reference{prl}.
\end{caselist}
\qed

We consider $T\cap X_i$. It is a collection of annuli and disks. 

\begin{slemma}\label{lm4.3}
We can achieve by isotopy and proper choice of $T$ that
$T\cap X_i$ is either empty, the whole $T$, or a single annulus.
Moreover, the intersection of an annulus $T\cap X_i$ with
$C$ is a pair of circles, each circle bounding a disk in $C\sm
T$ that contains exactly two of the 4 punctures $C\cap L$ of $P$.
\end{slemma}

For the proof let us fix a bit more language. Assume a torus
$T$ intersects a ball $X$ so that an annulus $A$ is a connected
component of $X\cap T$. Assume also the two circles in $\pa A$
are not contractible in $T$. (We will soon argue that this is
always the case.) One can only place two unlinked
unknotted not contractible loops on a torus, if they are two
meridians, or two longitudes and the torus is unknotted. Since
meridians (resp. longitudes) bound a disk only in the interior
(resp. exterior) of a solid torus, we can choose one (and only
one) of the complementary regions $Y$ of $T$ as the interior of
$T$ so that the loops $\pa A$ collapse in $Y$.

We then choose one of the two regions $Y'$ of $X\sm A$ so
that $Y'\cap C$ is a pair of disks (rather than an annulus).
By Sublemma \reference{sl40}, both disks intersect $L$ in
exactly 2 of the punctures each.
($T$ may enter into $Y'$, so that not necessarily $Y'=X\cap Y$.)
We call $Y'=\ir A$ the interior of $A$, and the exterior of $A$
is then obvious. Then $Y'$ is a cylinder. We call $A$ \em{(un)knotted}
if the core of $Y'$, or alternatively the intersection of a longitude
of $T$ with $A$, is a(n un)knotted arc in $X$. Similarly $T$ is
\em{(un)knotted (in $X$)} if $X\cap T=A$ and $A$ is (un)knotted. With
the same meaning we use this term when $X=X(Y)$ is a tangle space
and $A$ is disjoint from the tangle $Y$. (Then $\ir A\cap Y\ne\vn$
in general, and knottedness of an arc is understood as w.r.t. the
ball $B(Y)=X\cup Y$.) Note that $T$ is unknotted in a ball (but not
tangle space) $X$ if and only if $A$ is boundary parallel to $X$.

\mbox{}\hbox to \textwidth{%
\begin{minipage}{0.72\textwidth}
\parskip 5pt plus 3pt minus 2pt\relax

\parindent0pt\relax

We introduce a relation $\succ$ among annuli of the
considered type, saying for two such annuli $A,A'$ that
$A\succ A'$, if $A\subset \et A'$. It is easy to see that
this defines a partial order. (Beware, though, that this is
\em{not} equivalent to $\ir A\supset A'$, and this latter
condition is not reflexive.) A maximal element in $\succ$
is called an \em{outermost} annulus. 

Consider the example diagram on the right. It shows a view of $B(Y)$
from an equatorial section. The tangle $Y$ is depicted by the
thicker lines; the thinner lines indicate $C$ and $\pa T$.
The gray regions belong to $\ir T$. Then $A_1\succ A_0$ and
$A_3\succ A_2\succ A_0$, but $A_1$ does not compare to $A_{2,3}$.
However, $A_1\subset \ir A_2$ and also $A_2\subset \ir A_1$
(and the same is true for $A_3$ instead of $A_2$). The outermost
annuli are $A_{1,3}$.
\end{minipage}\hfill
\diag{1cm}{3}{3.6}{
  \pictranslate{1.5 1.8}{
    \picPSgraphics{2 setlinejoin}
    \picfilledcircle{0 0}{1.5}{}
    \picfillgraycol{0.8}
    \picfill{
      \picline{1.5 80 polar}{1.5 -80 polar}
      \picarcto{1.5 -100 polar}{-1.5}
      \piclineto{1.5 100 polar}
      \picarcto{1.5 80 polar}{-1.5}
    }
    \picmultigraphics[S]{2}{-1 1}{
      {\piclinewidth{14}
       \picline{0.12 -1.8}{0.12 1.8}
      }
      \picfill{
        \piccurve{1.5 30 polar}{1 0.5}{0.5 1}{1.5 60 polar}
        \picarcto{1.5 30 polar}{-1.5}
      }
      \picfill{
	\picarc{1.5 15 polar}{1.5 -65 polar}{1.5}
	\picarcto{1.5 -55 polar}{1.5}
	\picarcto{1.5 -0 polar}{-1.0}
	\picarcto{1.5 15 polar}{1.5}
      }
    }
    \picputtext{0.6 0}{$A_2$}
    \picputtext{-0.5 0.1}{$A_0$}
    \picputtext{0.7 0.7}{$A_1$}
    \picputtext{-0.7 0.7}{$A_1$}
    \picputtext{1.10 -0.62}{$A_3$}
  }
}
\hfill}
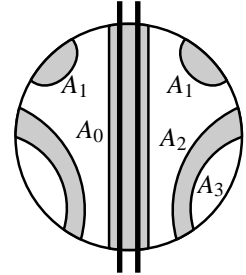
\vspace{1mm}

\proof[of Sublemma \reference{lm4.3}]
There is easily seen to be no separating disk of $Y_i$ in
$X_i$, so one can remove from $B_i$ all disks from $T\cap X_i$,
together with any other parts of $T$ in $B_i$ that lie on one
side of such disks. Then $T\cap X_i$ consists only of annuli.
(They are finitely many by compactness.)

If one of the circles in $T\cap C$ bounding an annulus $A$ of $T
\cap X_i$ is contractible in $T$, then $A$ is contained in a disk $D$
that is isotopable into the exterior of $T$ and not intersecting $L$.
Since $T$ is incompressible, the disk $D$, and hence $A$, is parallel
to $T$, and so $A$ can likewise be removed from $T\cap X_i$. So we can
assume that both circles in $\pa A$ are not contractible in $T$.
So we have the situation, and terminology available, discussed
before the proof.

Now we would like to rule out the possibility of several annuli
in $T\cap X_i$. For this assume w.l.o.g. that among all
essential tori $T$ of $L$, ours is chosen so that $T\cap
P$ has the fewest number of components (circles).

By the above argument, each annulus in $T \cap X_i$ bounds 
in $P$ a pair of meridional disks (with respect to one of
the complementary solid tori if $T$ is unknotted). In particular,
all the circles in $T\cap P$ are meridians of $\bd S=T$,
w.r.t. the interior $S=\ir T$ of $T$, or a proper choice of interior
if $T$ is unknotted. (Because a longitude and meridian always
intersect, the choice of $S$ cannot be different for different
circles in $T\cap P$.)

By Sublemma \reference{sl40}, each circle of $T\cap P=T\cap C$ which
bounds a disk in $C$ disjoint from $T\cap C$ (let us call such
circles \em{innermost}) intersects $\ge 2$ of the punctures
$L\cap C$ of $P$. There are clearly at least two circles in
$T\cap C$, and hence there are also at least two innermost.
Since $P$ has four punctures, we see that there must be
exactly two innermost circles, each bounding a disk in $C$
intersecting $L$ in exactly two punctures. Then $S\cap P$ is a 
collection of two twice-punctured disks, and unpunctured annuli.
Next we show that we can get disposed of the annuli in $S\cap P$.

Let $A$ be an annulus of $S\cap P$. Then $A$ forms a torus
$T_{1,2}$ with each of the two annuli that $\bd A$ cuts $T$ into. 
The $T_i$ inherit meridians from $T$, and their interior is defined
again as the region where meridians collapse. Then $\et T_i$
is determined also, $\et T=\et T_1\cup \et T_2$ and
$A=\et T_1\cap \et T_2$. We claim that at least one of
$T_{1,2}$ is essential. Since $A$ can be pushed into either $X_1$
or $X_2$, we have then a contradiction to the above minimizing
choice of $T$. 

First, $T_{1,2}$ do not compress in their interior, because $T$
does not. If some $T_j$ (is unknotted and) compresses in its
exterior, then all components of $L$ contained in $\et T_j$
lie within a ball contained in $\et T_j$. If there are such
components, $L$ is split, and otherwise, $T$ is isotopic to
$T_{3-j}$, and subsequently $A$ can be removed. 

If some $T_j$ were $\bd$-parallel to a component of $L$ in its interior
then $T$ would also be (and $T$ and $T_j$ would be isotopic).
Finally, at least one of $T_{1,2}$ is not $\bd$-parallel in
its exterior. If both were such, then because of $\et T_j
\subset \et T$, we would have both two components of $L$ in
the exterior of $T$, in contradiction to $S\cap L\ne\vn$.

With this argument we showed that any annulus in $S\cap C$
(that comes from a pair of nested annuli in $T\cap X_i$)
can be removed by isotopy. Thus we can achieve that $S\cap C$
consists only of disks. Also, by Sublemma \ref{sl40}, we argued
that there is only one pair of disks, so we have only one
annulus in $T\cap X_i$, and complete the proof of Sublemma
\ref{lm4.3}. \qed

We consider the two options for $T\cap C$ from Sublemma \ref{lm4.3}.

\begin{caselist}
\case \label{caseI}
$T\cap C\ne\vn$, that is, both $T\cap X_i$ are annuli.

\begin{slemma}\label{T1k}
$T$ is unknotted in $X_1$.
\end{slemma}

\proof Assume that $T$ is knotted in $X_1$. Then
$T\cap X_1$ is not parallel to the boundary of a string
of $Y_1$.  (Otherwise, the intersection of $T$ 
with $P$ is a pair of circles, each circle has only 
one, and not two as assumed, of the 4 punctures.)
If $d\ge 7$, then $U$ in \eqref{Tds} is not a rational tangle, and 
then $Y_1$ is not among the tangles in Theorem 4.9(a-d) of 
\cite{Wu}. This theorem says then that $T$ is simple, so excludes
such an annulus $T\cap X_1$. 

If $d=5$, then $Y_1$ is equivalent
(in the sense of definition \reference{dtg})
to a Montesinos tangle $M(1/2,p/q)$ with $q$ odd. To obtain a
contradiction in this case, assume w.l.o.g. $Y_1=M(1/2,p/q)$.
Let $Y_3$ be a prime tangle such that $L'=Y_3\cup Y_1$ is a
prime link of $\ge 2$ components.  Let $A \subset B(Y_3)$ be an
unknotted annulus identifying both circles of $T\cap P$ such
that it contains $Y_3$ in its interior. Consider the torus
$T'$ in $X(L')$ obtained by gluing $A$ and $T\cap X_1$.
So $T'$ is knotted. Let $S'$ be its interior. Then if $T'$ is 
$\bd $-parallel, it must be $\bd $-parallel to a single link 
component in $S'$. But since $L'$ has several components, and 
$S'$ contains all of $L'$, this is impossible. Since $T'$ is 
knotted, if it is compressible, then a compressing disk must
be meridional. Such a disk can be moved completely into either
$X_1$ or $X_2$, using that $Y_i$ have no separating disks.
But both is excluded, since $Y_{1,2}$ are prime and $P\cap L'$
is non-empty. Therefore, $T'$ is essential, and $L'$ is toroidal.

So any prime link $L'=Y_3\cup Y_1$ of $\ge 2$ components
is toroidal. To see that this is not so, take $Y_3=Y_1$. 
Since we do not know which pairs of punctures the two
circles of $T\cap P$ enclose, to glue the two annuli 
properly, we may need to rotate the two copies of $Y_1$
by $\pi/2$ or add a $\pm 1$ tangle. However,
in all cases these modifications can be
carried out so that $L'$ becomes a Montesinos link of 
length 4. (This observation will be required and
implicitly applied again in some of the below arguments.) 
Corollary 5 of \cite{Oertel} shows that such
links are atoroidal except if $p/q\ne \pm 1/2$,
which is clearly not the case here (because $q$ is odd). \qed

But now recall that $Y_1$ has no closed component. Then by
Sublemma \reference{lm4.3}, all of $Y_1$ lies in the
interior of $T$, i.e. in $S\cap B_1$. Since $T$ is unknotted,
it must be then $\pa$-parallel to $C$, and can be removed from
$X_1$. Thus it suffices to deal with the next case.

\case\label{caseII}
$T\cap C=\vn$. So $T$ lies in some $X_i$.
In our situation $Y_1$, $Y_2$ are, if not simple, up 
to equivalence $M(1/2,p/q)$. So it suffices that we study the
case $Y_1=M(1/2,p/q)$ (with $q$ even or odd, that is, with or 
without closed component) and assume $T\subset X_1$.

We obtain by inclusion a torus $T$ in the exterior of the
link $L'=Y_1\cup Y_3$ for any tangle $Y_3$. Again we want
to obtain a contradiction from this by choosing $Y_3$ well
and using Oertel. Assume $Y_3$ is prime and $L'$ is non-split.

We claim that this torus $T\subset X_1$ is not compressible
in $E(L')$. To see this, assume $T$ were compressible. First note
that if $T$ separates components of $Y_1$ in $X_1$, it would too
in $L'$, in contradiction to the non-splitness of $L'$. So 
$T$ separates no components in $X_1$. Then the only way in which
$T$ would be compressible in $E(L')$ but incompressible in
$E(L)$ is that $T$ is knotted, and $X_1\supset \et T$.

Let $D$ be a compressing disk of $T$ in $E(L')$. This disk
may penetrate into $X_3=X(Y_3)$. But since $Y_3$ was
chosen prime, $X_3$ has no separating disks, and so $D$ can
be moved out of $X_3$, and into $X_1$. So $T$ would
compress in $X_1$ too, a contradiction.

With this we assure that $T\subset X_1$ is incompressible in
$E(L')$. So it is essential, unless it is boundary parallel. 
It is not boundary parallel to a closed component of $Y_1$,
because it is essential in $E(L)$, and so also in $X_1$. So 
$T$ can only be boundary parallel in its region containing
$B_3$. This can be avoided for example by choosing $Y_3$
to have a closed component. 

Therefore, all $L'=Y_1\cup Y_3$ where $Y_3$ has a closed
component must be toroidal. This is easily disproved by choosing
$Y_3$ well (so that $L'$ is a Montesinos link) and using Oertel.
\end{caselist}

Since we obtained contradictions in all cases, we conclude that
there is no $T$, and Lemma \reference{lato} is proved. \qed

%


\section{Links of more components\label{S51}}


Now we derive some consequences and generalizations for links
of more components. (In \S\ref{S51}
we use consistently $n=n(L)$ for the number of components of
a link $L$ and $g=g(L)$ for its genus. The cases $n(L)\le 2$
were discussed before, so assume throughout $n\ge 3$.) 

The first theorem deals with fiberedness.
Kanenobu \cite{Kanenobu} extended the realization of monic
polynomials to fibered links. However, his construction, which
seems the only one known, uses connected sum with Hopf links.
Thus, for more than two components, surprisingly, the simple
question to find a prime fibered link appears open (even for $n=2$,
Kanenobu's links are not proved to be prime). The theorem removes
this shortcoming, with a more specific statement.

\begin{theorem}\label{HHHH}
Let $\nb$ be an admissible (as in definition \ref{dad})
monic Conway polynomial of an $n$-component link, $n\ge 3$.
Then, except for $n=3$, $g(L)=0$ and $\nb=+z^2$, there exists
a prime arborescent fibered link $L$ with $\nb_L=\nb$, such that
the fiber of $L$ is a canonical surface obtained from a special
arborescent diagram of $L$. Unless $n=3$, $g(L)=0$, and
$\nb=-z^2$, the link $L$ is hyperbolic, and
\[
\vol(L)\,\le\,10V_0\cdot\,\left\{\,\begin{array}{c@{\quad}l}
2\Md\nb-n & \mbox{if }g(L)>0 \\[1mm] n & \mbox{if }g(L)=0
\end{array}\,\right..
\]
\end{theorem}

The following object will be useful for the primeness and
hyperbolicity arguments.

\begin{defi}
Define the \em{linking graph} $G(L)$ of $L$ by putting a
vertex for each component of $L$ and connecting vertices of
components with non-zero linking number. Optionally, we
may label an edge by the linking number.
\end{defi}

\proof[of Theorem \reference{HHHH}] Let first $g(L)>0$.
We deal with the case $n=3$ first. Consider the 2 component 
link $L'$ found in Theorem \ref{th2c} for $\nb'=\nb/z+z$.
Assume the (reverse) clasp ** in the left diagram of figure
\ref{figgr} is negative, by possibly mirroring $L'$ (mirroring
preserves $\nb$ for even number of components). Recall that
$L'$ is obtained from a knot as on the left of figure \ref{figgr}
by smoothing out one crossing in its parallel clasp *.

Call the replacement of a crossing with a parallel clasp a
\em{clasping}, and give it a sign as for the crossings involved:
\begin{eqn}\label{clasping}
\diag{7mm}{1}{1}{
    \picmultivecline{0.18 1 -1.0 0}{0 1}{1 0}
    \picmultivecline{0.18 1 -1.0 0}{0 0}{1 1}
}\quad\lra\quad
\diag{7mm}{2}{2}{
  \picPSgraphics{0 setlinecap}
  \pictranslate{1 1}{
    \picrotate{-90}{
      \rbraid{0 -0.5}{1 1}
      \rbraid{0 0.5}{1 1}
      \pictranslate{-0.5 0}{
      \picvecline{0.02 .95}{0 1}
      \picvecline{0.98 .95}{1 1}
    }
    }
    }
}\es,\rx{1.5cm}
\diag{7mm}{1}{1}{
    \picmultivecline{0.18 1 -1.0 0}{0 0}{1 1}
    \picmultivecline{0.18 1 -1.0 0}{0 1}{1 0}
}\quad\lra\quad
\diag{7mm}{2}{2}{
  \picPSgraphics{0 setlinecap}
  \pictranslate{1 1}{
    \picrotate{-90}{
      \lbraid{0 -0.5}{1 1}
      \lbraid{0 0.5}{1 1}
      \pictranslate{-0.5 0}{
      \picvecline{0.02 .95}{0 1}
      \picvecline{0.98 .95}{1 1}
    }
    }
    }
}
\es.
\end{eqn}

Then apply a positive clasping at a crossing among those
corresponding to the edge labeled $2a_2-1$ in figure \ref{figgr}.
(If these crossings are negative, create a trivial clasp by a
Reidemeister II move in advance.) We claim that the resulting
3-component link $L$ is what we sought.

The Conway polynomial is easily checked using the skein relation
\eqref{sRs} at the crossing created by the clasping. In that case
$D_+$ depicts $L$, $D_-$ depicts a
$(2,-2,k)$-pretzel link ($k$ even), and $D_0$ depicts $L'$. By
the proper choice of $\nb'$, we see that $\nb_L=\nb$.

The fibering is also easy, since a clasping results in a Hopf
plumbing on the canonical surface. By \cite{Gabai2,Gabai3},
the fiber property is invariant under a Hopf plumbing.

It remains to see primeness. This can be shown again from the
arborescency using \cite{KL}, but there is a more elementary
argument. Note that all components of $L$ are unknotted and have
pairwise non-zero linking number. (Here the proper choice of
signs of clasps is helpful.) Thus if we had $L=L_1\# L_2$, 
former property excludes the option that some of
$L_{1,2}$ is a knot, and latter property excludes the
option that both are 2-component links.

For $n>3$ we can use induction. Again we apply
claspings (either sign may do) at some of the crossings of
$2a_2-1$ (possibly creating new crossings by Reidemeister
II moves). The link on the bottom right of figure \ref{figgr}
is a typical example (for $n=4$). Again the check of $\nb$
is easy; $D_\pm$ depicts $L$, $D_0$ depicts a connected
sum of a $(2,-2,k)$-pretzel link with Hopf links, and $D_\mp$
depicts a link of the sort constructed for $n-1$.
The skein relation of $\nb$ again easily allows
to adjust the polynomial of $L$ properly.

To see primeness, use again that all components of $L$ are
unknotted. So if $L=L_1\# L_2$, then both of $L_{1,2}$ are
links. Then the linking graph $G(L)$ of $L$ must
have a cut vertex $v$ (i.e. it must become disconnected when
removing $v$ and its incident edges). However, for our $L$ this 
is easily seen not to be the case. Here $G(L)$ consists of a
chain connecting all vertices, with an additional edge between
two vertices of distance 2 in the chain. So $L$ is prime. Let
us display the graphs for $3$ and $4$ components, also for
future reference. They look like (up to reversing sign in
all linking numbers)
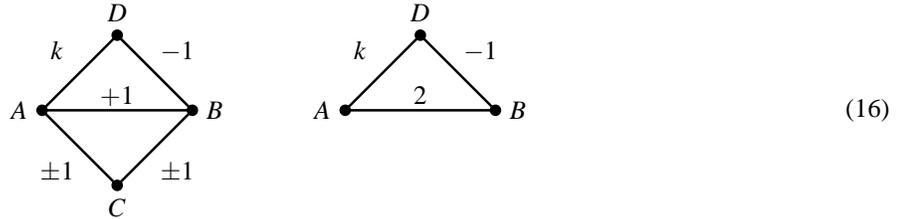
\begin{eqn}\label{lg}
\diag{1.0cm}{2}{2.8}{
 \pictranslate{0 0.4}{
  \tycl{1 2}{2 1}{0 1}{1 0}
  \picputtext{1 1.2}{$+1$}
  \picputtext{1.8 1.8}{$-1$}
  \picputtext{1.8 0.2}{$\pm 1$}
  \picputtext{0.2 0.2}{$\pm 1$}
  \picputtext{0.2 1.8}{$k$}
  \picputtext{-0.3 1}{$A$}
  \picputtext{2.3 1}{$B$}
  \picputtext{-0.3 1 x}{$C$}
  \picputtext{2.3 1 x}{$D$}
 }
}\kern2cm
\diag{1.0cm}{2}{2.8}{
 \pictranslate{0 0.4}{
  \trig{2 1}{0 1}{1 2}
  \picputtext{1 1.2}{$2$}
  \picputtext{1.8 1.8}{$-1$}
  \picputtext{0.2 1.8}{$k$}
  \picputtext{-0.3 1}{$A$}
  \picputtext{2.3 1}{$B$}
  \picputtext{2.3 1 x}{$D$}
 }
}
\end{eqn}
Here the component designation for $n=4$ is as in figure
\reference{figgr}.
Note that, since the diagram is special, $D$ and $A$
have with $B$ a linking of opposite sign.

For $n=3$ we let $C$ identify with $A$ under undoing
one of the claspings \eqref{clasping} in the $n=4$ case.
As occurred in the primeness argument, we can have also
$lk(A,B)=0$. We will need this case only once (at the
end of the proof of Lemma \reference{lXYz}), and otherwise
stick with $lk(A,B)=2$.

Our construction yields links with all desired properties
(except hyperbolicity, which we treat below)
whenever $g(L)>0$. Finally, turn to the
case $g(L)=0$. We use the pretzel links of type I in Gabai's
theorem 6.7 in \cite{Gabai4}. The links in case 1 (B), (C)
there realize the stated polynomials. For even number of
components, case (C) applies, and we get both possible polynomials
$\pm z^{n-1}$ by mirroring (which changes sign of $\nb$). For
odd number of components we have the pretzel links in case (B).
To see that their polynomial is $(-1)^{\br{n/2}}z^{n-1}$, one
can use, for example, the formula of Hosokawa-Hoste \cite{Hoste}.
For $n=5,7,\dots$ and $\nb=(-1)^{\BR{n/2}}z^{n-1}$ we found,
with the help of some computation, the sequence of links
with Conway notation $(2,2,-2)(2,-2,2,\dots,-2,2)$, the first two
of which look like:
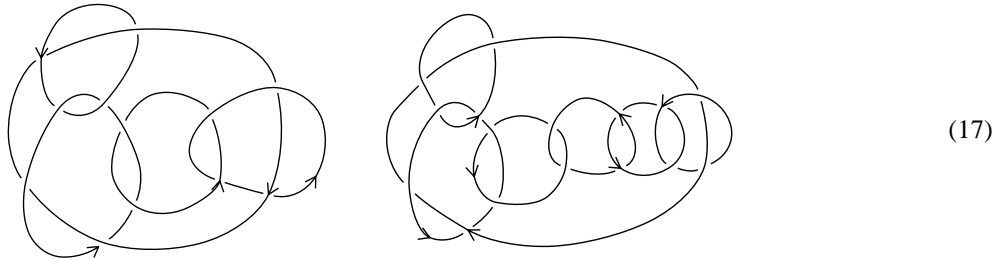
\begin{eqn}\label{uu}
\epsfsv{3.6cm}{t1-16_409274_lnk_nice}\qquad
\epsfsv{3.6cm}{t1-16_409274_lnk2_nice}
\end{eqn}
(The orientation of components is so that all clasps are
reverse.) The fibering of these examples can be confirmed by the
disk (product) decomposition of Gabai \cite{Gabai4}, and the
proper $\nb$ using \cite{Hoste}.

We postpone the hyperbolicity proof to lemmas \ref{lXYz} and \ref{lXY}.
The volume estimate is again easy from Theorem \reference{LAT}. \qed

\begin{rem}\label{R51} The following observations indicate how
one can (or can not) modify or extend Theorem \ref{HHHH}.
\begin{mylist}{\arabic}
\myitem For $n=3$ the only diagrams with canonical surfaces of genus
$0$ are the $(p,q,r)$-pretzel diagrams, $p,q,r$ even. Then Theorem
6.7 Case (1) of Gabai \cite{Gabai4} shows that there is no prime
link for $\nb=+z^2$, even with a canonical fiber surface from an
arbitrary diagram. 

\label{L51}
\myitem 
The algebraic topologist considers $\Dl$ usually up to units
in $\bZ[t^{\pm 1}]$, in opposition to treating $\Dl$ as the
equivalent \eqref{nbDl} of $\nb$. In that weaker sense the
exceptional links \eqref{uu}
in our proof could be avoided. For knots the ambiguity of $\Dl$
is not essential, because the condition $\Dl(1)=1$ allows one
to recover the stricter form. Note, though, that for links of
more than one component, we lose the information of a sign
in the up-to-units version.

\myitem\label{i4,}
The exception $n=3,g=0$ also disappears for the strict $\Dl\ne 0$
if we waive on fiberedness (and then also on monic polynomials)
and demand $2\Md\Dl=1-\chi$ instead.
The corresponding statement follows just by an obvious
modification of the proof we gave. (For genus 0 one can easily
adjust infinitely many pretzel links to give the proper polynomial.)

\end{mylist}
\end{rem}

If we like to keep small 4-genus, we have

\begin{corr}\label{cor1}
For any admissible Alexander polynomial $\Dl$ of a link, there exists
an arborescent link $L$ with $\Dl(L)=\Dl$, $\Md\Dl=1-\chi_c(L)$ and
$\chi_s(L)\ge -1$. Moreover, $L$ can be chosen to be fibered
if $\Dl$ is monic.
\end{corr}

\begin{rem}
Clearly for an $n$-component link, $\chi_s\le n$,
but even below this bound, one cannot augment
$\chi_s$ unrestrictedly, since it is related to
(the vanishing of) certain linking numbers, which
in turn have impact on the low-degree terms in $\nb$.
(In particular $\chi_s=n$ means strongly slice,
which implies that $\nb=0$.)
\end{rem}

\proof For one component, $u(K)\le 1$ implies $\chi_s(K)\ge -1$.
For a link of two components take the link constructed
for Theorem \reference{th2c}. Observe that this link bounds a ribbon
annulus, so $\chi_s\ge 0$. For $n\ge 3$ components,
we can always achieve that $\chi_s\ge -1$ for the links $L$
in Theorem \reference{HHHH}, by varying the sign of claspings
\eqref{clasping} with the parity of $n$. \qed


\begin{lemma}\label{lXYz}
The link $L$ of Theorem \ref{HHHH} is choosable to have a complement
which is not Seifert fibered, unless $n=3$, $g(L)=0$, and $\nb=-z^2$.
\end{lemma}

\proof  Consider first $g(L)>0$.
We use again the description in \cite{BM}.
Since $n\ge 3$, all components are unknotted, we have only the
types shown in figures 2 (type (a)) and 3 (type (b)) therein.
Now all these links have the following property: there is
a component $M$ having the same linking number with all the
others, up to sign. Looking at $G(L)$ for our links $L$, we see
that only $n\le 4$ components come in question. 

So for $n=4$, $M$ can be only one of $A$ or $B$ (see \eqref{lg}).
However, the next property of Burde-Murasugi's links is that
all components different from $M$ have mutually the same linking
number. This immediately rules out also $n=4$.

Now for $n=3$, $M$ can be only $D$ and $k=\pm 1$. In
type (b) of Burde-Murasugi, the distinguished component
$M$ has linking number $\pm \ap$ with all the other
components, and in that case it was assumed that $\ap>1$,
so this option is ruled out. It remains their type (a).
For these links, looking at Figure 2 of \cite{BM} with
$m=3$, and taking care of linking numbers, we see
that we have the $(2,-2,4)$-pretzel link, oriented
so as to be the closure of the 3-braids
$\sg_2^{-1}\sg_1^{-2}\sg_2^{-1}\sg_1^{\pm 4}$, but
for $\sg_1^{-4}$ one component involving these
crossings must be reversed. Latter case gives a link
of genus 0, so consider only former, i.e. with $\sg_1^4$
in the braid.

The Conway polynomial of this link is $\nb=-3z^2-z^4$. The link $L$
(up to mirroring) obtained from our construction with such polynomial
is shown on the left of \eqref{yu}. It has the linking graph on
the right of \eqref{lg} for $k=-1$. It turns out that SnapPea reports
this link non-hyperbolic, so apparently it is the Burde-Murasugi link.
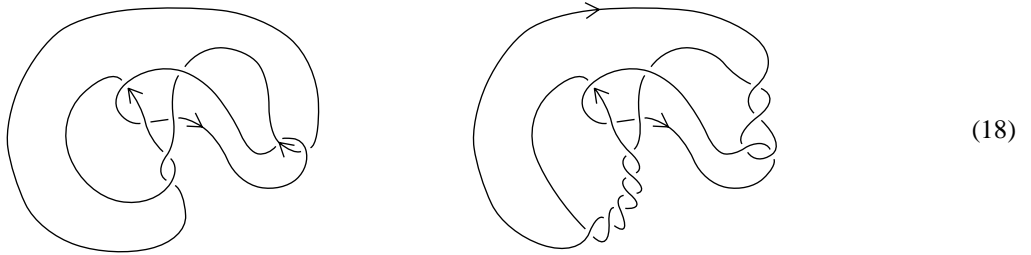
\begin{eqn}\label{yu}
\epsfsv{3.5cm}{t1-aparb_hyp90}
\kern2cm
\epsfsv{3.5cm}{t1-aparb_hyp90a}
\end{eqn}
However, now recall that we had some option in the construction
of $L$. First we can change the sign of the clasp * in \eqref{gj},
which here leads to a composite link. Next, though, we can
change the sign of the clasping \eqref{clasping}. This leads
to another link with the same polynomial, given on the right
of \eqref{yu}. SnapPea reports it to be hyperbolic, with which
the case $g(L)>0$ is finished.

The links $L$ of genus $0$ are dealt with by the same
argument. Again by linking numbers we are down to $4$
components (in particular all those links of \eqref{uu}
are done). For 4 components, the linking graph of a pretzel is
a cycle of length 4, so this case is out too, and for $n=3$
we arrive at the additional exception we had to make~-- the
$(2,-2,r)$-pretzel links are indeed Seifert fibered. \qed

\begin{lemma}\label{lXY}
The link $L$ of Theorem \reference{HHHH} is atoroidal.
\end{lemma}

\proof[of Lemma \reference{lXY}] Let us focus on $g(L)>0$.
We adapt the proof, as far as possible, from lemma
\reference{lato}, and use the notation from there.
The tangle decomposition of $L$ in \eqref{Tds} modifies
so that now
\begin{eqn}\label{Tds'}
Y_1\ =(U\ 1\ 1,-2)\,m\quad\mbox{and}\quad Y_2\ =\ ((2k,-2)\ 1\ 1\ ,
\pm 2,\pm 2,\dots,\pm 2)
\end{eqn}
Again let $T$ be an essential torus in $E(L)$.
Since both $Y_{1,2}$ are again easily proved to
be prime, we can assume w.l.o.g. that $T$ does not intersect
any tangle space $X_i$ in disks, but only in annuli. Still
$T_1$ has no closed component and is subjectable to \cite{Wu}.
Then all intersections of $T$ with the tangle sphere $C$ are
meridional disks, with respect to a proper choice of interior
$S=\ir T$. Assume again $T$ is chosen so that $S\cap C$ has
the fewest connected components.

With the same argument we have first:

\begin{slemma}\label{sl4.1'}
Sublemma \reference{sl40} holds. \qed
\end{slemma}

\begin{slemma}\label{sl4.2'}
If $T$ is knotted, then $T$ is not separating, i.e. $L\subset
\ir T$. 
\end{slemma}

\proof All components of $L$ are unknotted. Any unknot embedded
in a knotted solid torus has homological degree $0$. So for
each pair of components $M_1\in\ir T$, $M_2\in\et T$, we must
have $lk(M_1,M_2)=0$. So if $T$ were separating, the linking
graph $G(L)$ would be disconnected, which we saw is not the
case. Since by incompressibility, there is always some $M_1$,
there cannot be any $M_2$. \qed

\begin{slemma}\label{sl4.3'}
If $T$ is unknotted, then $T$ is separating. Let $P,Q$ be the
sets of components of $L$ in $\ir T$ resp. $\et T$. Then
$G=G(L)$ has the following property. If for some $a\in P$,
$b\in Q$ there is no edge between $a$ and $b$ in $G$, then
there is no edge between $a$ and $b'$ for any $b'\in Q$, or
there is no edge between $a'$ and $b$ for any $a'\in P$.
\end{slemma}

\proof Clearly an unknotted torus must separate, else it
would compress. Now when $T$ is unknotted, $L$ is a satellite
of the Hopf link. Then for two components $a\in P$, $b\in Q$
of $L$ we have $lk(a,b)=[a]\cdot [b]$, where the brackets
denote the homology class in $H_1(\ir T)=H_1(\et T)=\bZ$.
So if $a$ and $b$ are not connected in $G$, one of $[a]$ or
$[b]$ must be $0$, and the claim is clear. \qed

\begin{slemma}\label{sl4.4'}
Sublemma \reference{lm4.3} holds still.
\end{slemma}

\proof The proof of Sublemma \reference{lm4.3} goes through
with the help of now Sublemma \reference{sl4.1'}, except for the
argument why some of $T_{1,2}$ is not $\pa$-parallel in its
exterior. 

If $T$ is knotted, then by Sublemma \reference{sl4.2'},
its exterior is empty, so clearly none of $T_{1,2}$ can be
$\pa$-parallel in its exterior. If $T$ is unknotted, then all
annuli of $T\cap X_i$ are unknotted too. 
Now since one of the $Y_i$, namely $Y_1$, still has
no closed component, an outermost annulus of $T\cap X_1$ is parallel
to $C$. Then successively all annuli of $T\cap X_1$ can be removed,
so $T\cap X_1=\vn$. \qed

Back to the proof of lemma \reference{lXY},
now we can apply \cite{Wu} to $Y_1$. An annulus $T\cap X_1$
must be parallel to $C$, provided $U$ in \eqref{Tds'} is
not a rational tangle. Then $T$ can be removed from $X_1$, so
$T\subset X_2$. If $U$ is rational and $T\cap X_1\ne\vn$, then
we can obtain a contradiction to Oertel's result by joining
$Y_2$ with itself properly to obtain a Montesinos link of
length $4$.

So we can assume $T\subset X_2$.

%
%
%

Now let $L'=Y_3\cup Y_2$ be a prime (non-split) link of $\ge 5$
components, and $Y_3$ be a prime tangle with a closed component.
We claim that $T\subset E(L')$ is essential. The argument is the
same as in case \ref{caseII} of the proof of Lemma \ref{lato}.
So again all such $L'$ would be non-atoroidal.

Thus we can conclude the proof of Lemma \reference{lXY} for $g>0$
with Lemma \reference{lXY'} below. For our links of $g=0$, we
can apply Oertel to the pretzel links,
and the links in \eqref{uu} are dealt with
the same argument is those in Lemma \reference{lXY'}. (See
the remark at the end of its proof.) \qed

\begin{lemma}\label{lXY'}
The links $L'$ with Conway notation 
\[
((k,-2)\ 1\ 1),\pm 2,\pm 2, \dots,\pm 2,0\ m\,,
\]
of $n(L')\ge 5$ components for $k,m\in \bZ$, $k\ne 0$
even, are atoroidal.
\end{lemma}

Here is an example $L'$ with $m=0$, $k=4$ and $n=5$ components,
together with its linking graph $G(L')$ we will use shortly.
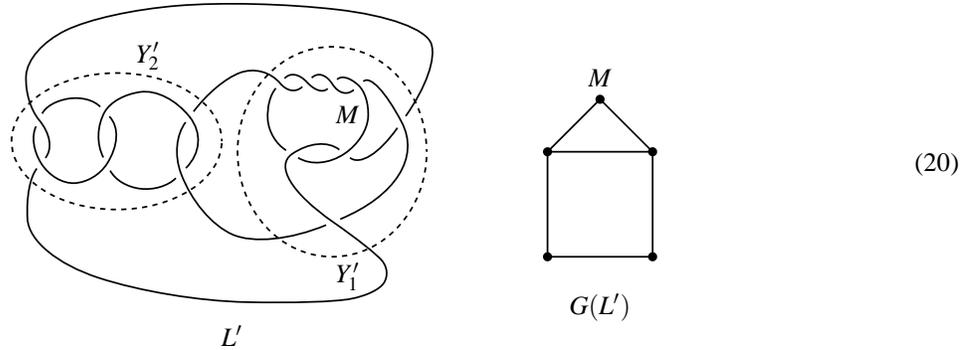
\begin{eqn}\label{lXY'x}
\diag{7mm}{12}{7.5}{
  \picputtext{4 4}{
    \epsfs{4.6cm}{t1-aparb_hyp91}
  }
  \picputtext{6.2 4.7}{$M$}
  \picputtext{4 0.5}{$L'$}
  {\piclinedash{0.1}{0.05}
   \picellipse{1.8 4.2}{2 1.3}{}
   \picellipse{5.9 4}{1.8 2}{}
  }
  \picputtext{2.4 5.9}{$Y_2'$}
  \picputtext{6.2 1.7}{$Y_1'$}
  \pictranslate{1 -1}{
    \picputtext{10 6.4}{$M$}
    \picputtext{10 2.1}{$G(L')$}
    \ttycl{10 6}{9 5}{11 5}{9 3}{11 3}
  }
}
\end{eqn}

\proof Let $Y_1'=(k,-2)\ 1\ 1\ m$ and $Y_2'=(\pm 2,\pm 2,\dots, \pm
2)$. Then $L'=Y_1'\cup Y_2'$. (Follow the diagrams in \eqref{lXY'x}.)

If we remove the closed component $M$ of $Y_1'$, then we have a
pretzel link, which is atoroidal by Oertel. (Here we may better avoid
the $(2,-2,2,-2)$-pretzel link $L'\sm M$; but we will just see that
its unique essential torus still fits into the below conclusions.) Thus
an essential torus $T$ of $L'$ must become inessential in $L'\sm M$.
Since $L'$ is non-split, this means that one of the regions of $T$
must contain either only $M$ (if $T$ compresses in $L'\sm M$), or
$M$ and exactly one other component $M'$ of $L'$ (if $T$ is
$\pa$-parallel to $M'$ in $L'\sm M$). In particular, since we
have $n\ge 3$ components, $T$ is separating.

\def\hf{\myfrac{1}{2}}
Now again all components of $L'$ are unknotted and $G(L')$
is connected. So $T$ separating means by Sublemma \ref{sl4.2'}
that $T$ is unknotted (as for the essential torus of
$M(\hf,-\hf,\hf,-\hf)$). Now we can apply Sublemma \ref{sl4.3'}
on $G(L')$. For $n(L')\ge 5$ components, we easily see that
the option $T$ containing a component $M'\ne M$ is ruled out. 

Thus $T$ contains $M$ alone in one region (and $n-1\ge 4$
components of $L$ in the other one). Then by the argument for
Sublemma \ref{sl4.4'}, $T$ can be isotoped (or chosen more properly)
into $X_1'=X(Y_1')$ or $X_2'=X(Y_2')$. Let us explain this briefly.

First, the argument excluding $T_{1,2}$ being both $\pa$-parallel
in their exterior applies now, because we assured that none
of the regions of $T$ contains precisely 2 components of $L'$.
So the conclusion of Sublemma \ref{lm4.3} applies. Next, the option
of an annular intersection $T\cap X_i'$ is ruled out as follows.

The annuli $T\cap X_1'$ and $T\cap X_2'$ again determine an
interior of $T$ by letting the circles in $T\cap C$ collapse
therein. Now $T$ is unknotted and contains only one component
in its exterior, a component which does not intersect $C$.
Then for at least one $i=1,2$ the annulus $T\cap X_i'$ will be
(unknotted and) with empty exterior in $X_i$, so $\pa$-parallel
to $C$, and could be removed. 

Now having $T$ within $X_1'$ or $X_2'$, we can obtain the same
contradiction as before by looking at $Y_1\cup Y_3$ or
$Y_2\cup Y_3$ for proper $Y_3$ and applying Oertel.

Let us say a word on the links in \eqref{uu}. Their
linking graph is the same as for our $L'$. Again
removing $M$, when specifying it so as the labeling
in graph on the right of \eqref{lXY'x} to be correct,
gives a pretzel link. So the argument here applies
unchangedly. \qed

%
Let us conclude the hyperbolicity proof with a few general/historic remarks.
One reason for the effort we needed to spend we see in the lack of
extension of Wu's work \cite{Wu} to tangles with closed components.
This extension is a substantial program, and we were forced to go
some steps along it, even though it was not our primary focus. It is
clear that our method can be applied to many more examples, although
the complete treatment of arborescent tangles is still far ahead.

The other main motivation for our hyperbolicity proofs was the status of
Bonahon-Siebenmann's monograph \cite{BS}. We were aware that we reprove
their theorem on the classification of hyperbolic arborescent links
in particular special cases. Still we were bothered by the notorious
inavailability of \cite{BS}, announced decades ago, but never
completed. Even for Montesinos links, written accounts needed some
amendment. At least atoroidality of the link complements seemed not
completely clarified. An additional complication for links is
that not only torus links have Seifert fibered complements.
Among the links in \cite{BM}, at least the $(2, -2, r)$ pretzel
links, pointed out by Ying-Qing Wu, are Montesinos and (for
$|r|\ne 1,2$) non-torus links whose complements are Seifert
fibered (and atoroidal). Thus in particular the statement and
proof of corollary 5 in \cite{Oertel} must be corrected
accordingly (see e.g. also \cite{adeq}).

Only after we completed our work, we were informed of a recent preprint
of Futer and Gu\'eritaud \cite{FG}, which gives a written proof of
Bonahon-Siebenmann's theorem characterizing the hyperbolic arborescent links.
Still it seems fair to say that our effort was (almost) simultaneous,
independent, shorter than the (full extent of the) work in
\cite{FG}, and makes our paper more self-contained. Thus we see both
some right and some sense to keep the material in \S\ref{S4} and
\reference{S51}, rather than mostly avoid it by referring to \cite{FG}.


\section{Tangle surgery constructions}

The following constructions, which are also heavily used in
\cite{metabol}, show infinite families of links with given
polynomial, if we focus on arborescency and $\chi_s$, but abandon
fibering and, in certain cases, minimality of the canonical surface.
(Note that, in \cite{canon} we showed that almost every monic Alexander
knot polynomial of degree 2 is realized by only finitely many
canonical fiber surfaces, so abandoning fibering of the canonical
surface is a non-trivial relaxation. See \S\reference{S72} for
related discussion.)

We will use some tangle surgery arguments. With the terminology
of Definition \reference{dtw}, we state first

\begin{lemma}\label{lsurg}
Let $S_k$, for $k\in\bZ$, $k\ne 0$, be the $(1,2k-1)$ pretzel
tangle, with orientation chosen so that the twist of $2k-1$ is
reverse. ($S_1$ is a positive parallel clasp.) Then $S_k$ can
be replaced by tangles $T_{p,q,r}$, that contain three twists
of $p,q,r$, such that all lengths $|p|,|q|,|r|$ can be chosen
arbitrarily large, and any such tangle replacement preserves the
Alexander polynomial.
\end{lemma}

\proof
%
Consider the $(p\pm 1,q,r)$-pretzel knot diagrams $D(p\pm 1,q,r)$,
with $p\pm 1,q,r$ odd. Their Alexander polynomial is determined by 
$v_2=\myfrac{1}{2}\Dl''(1)$, which is
\[
v_{2,\pm}\,=\,\frac{(p\pm 1)q+(p\pm 1)r+qr+1}{4}\,.
\]
Now for $p=0$, $q=1$, $r=2k-1$ we have 
\begin{eqn}\label{str}
v_{2,+}\,=\,k\,,\qquad v_{2,-}\,=\,0\,.
\end{eqn}
We need to find more solutions to \eqref{str}. We have
\begin{eqnarray}
(p-1)q+(p-1)r+qr+1 & = & 0 \label{(A)} \\
(p+1)q+(p+1)r+qr+1 & = & 4k \label{(B)} 
\end{eqnarray}
Then $\eqref{(A)}-\eqref{(B)}$ gives $q+r=2k$, and 
$\eqref{(A)}+\eqref{(B)}$ gives $p(2q+2r)+2qr=4k-2$,
so
\[
p=\frac{2k-1-qr}{2k}\,.
\]
We would like $p\in \bZ$ and $p$ even. To achieve this, choose
\begin{eqn}\label{xxx}
q\,=\,1+2nk\,,\qquad r=2k-1-2nk\,,
\end{eqn}
for $n\in \bZ$. Let $T_{p-1,q,r}$ be the tangle obtained by cutting
out from $D(p\pm 1,q,r)$ the switched crossing, for example 
for $(p,q,r)=(8,5,-3)$:
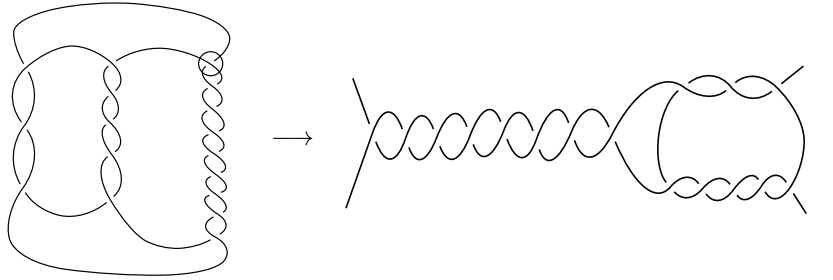
\[
\diag{8mm}{4}{5}{
  \picputtext{2 2.5}{\epsfs{3.5cm}{t1-P359}}
  \piclinewidth{50}
  \piccircle{3.55 3.74}{0.2}{}
}\quad\lra\quad
\epsfsv{3.6cm}{t1-dlsurg0}\quad.
\]
(The shift to make
the first index odd is done for future convenience.) Now we can
substitute $T_{p-1,q,r}$ for $S_k$, so that $\Dl$ is preserved (see
\cite{Bleiler} or \cite{STW}). Also $|p|,|q|,|r|\to\infty$ when
$|n|\to \infty$. \qed

\begin{rem}
We will use also the surgery on the mirrored tangles. The mirrored
surgery for $k=1$ and $(p,q,r)=(8,5,-3)$ is shown below:
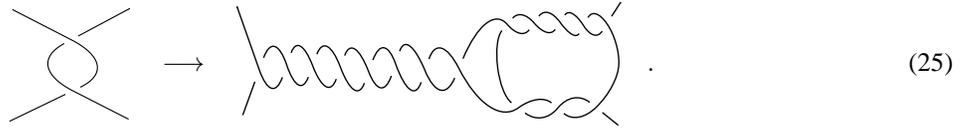
\begin{eqn}\label{tgs}
\epsfsv{1.6cm}{t1-dlsurg5}\quad\lra\quad
\epsfsv{3cm}{t1-dlsurg1}\quad.
\end{eqn}
\end{rem}

If we abandon fiberedness and relax the minimal genus condition
$2\Md\Dl=1-\chi$, then for example, we easily restore arborescency
in corollary \reference{cor1}:

\begin{corr}\label{cor1x}
For any admissible Alexander polynomial $\Dl$ of a link, there
exists an arborescent link $L$ with $\Dl(L)=\Dl$, $\Md\Dl\ge
-3-\chi_c(L)$ and $\chi_s(L)\ge -1$.
\end{corr}

\proof Consider the link in the proof of Corollary \ref{cor1}.
Let $D$ be the diagram constructed there. We apply the
modifications in \eqref{tsurg}. Create a prime diagram $D'$ by
adding a positive and negative parallel clasp. Then apply tangle
surgery on these clasps in $D'$ with two mutually mirrored tangles,
so that one obtains a diagram $D''$ of a concordant link. For
$(p,q,r)=(8,5,-3)$ the operation looks as follows:
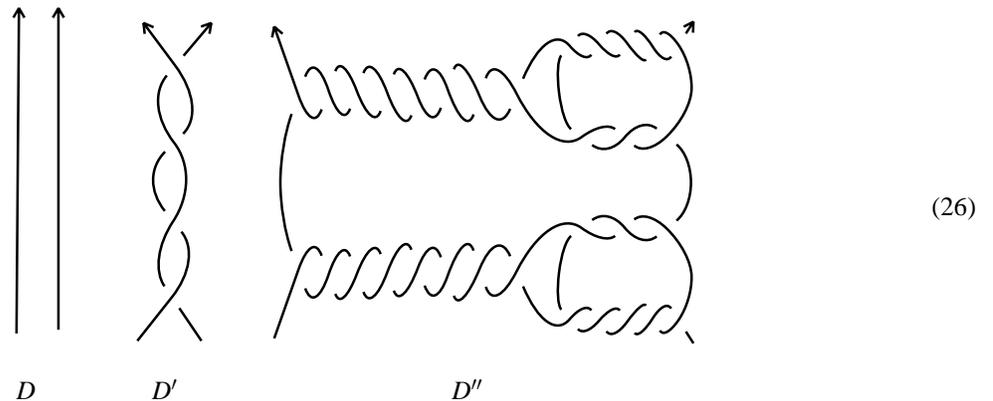
\begin{eqn}\label{tsurg}
\begin{array}{ccc}
\multicolumn{3}{c}{\epsfs{6cm}{t1-dlsurg}} \\
\ry{6mm}D\rx{7mm} & \rx{5mm}D' & \rx{5mm}D'' \\
\end{array}
\end{eqn}
These two tangle surgeries preserve arborescency and
$\chi_s$ and augment the genus of the diagram at most
by two. \qed

If we are interested in controlling only $\chi_s$, there are
virtually no difficulties at all in using surgeries, and we have:

\begin{corr}
For any admissible Alexander polynomial $\Dl$ of an $n$-component
link and $\chi\le 0$ with $n+\chi$ even, there exists an
arborescent link $L$ with $\Dl(L)=\Dl$ and $\chi_s(L)=\chi$.
\end{corr}

\proof The largest $\chi$ was dealt with in corollary \ref{cor1x}.
Then take iterated connected sum with $(-3,5,7)$-pretzel knots
and apply the (concordance) surgery \eqref{tsurg}. \qed

%
%



\section{Infinite families of links\label{S72}}

It is a natural question which admissible monic Alexander
polynomials are realized by \em{infinitely many} fibered
links. For knots the problem was suggested by Neuwirth
and solved fully by Morton \cite{Morton} (after previous
partial results; see for example Quach \cite{Quach}). As well
known, genus one fibered knots are only the trefoil and
figure-8 knot. In contrast, Morton constructs for each possible
monic Alexander polynomial of maximal degree greater than one an
infinite sequence of distinct fibered knots with this polynomial
(though without regard to any additional knot properties).

Unfortunately, extensions of Morton's construction to links
seem never to have been attempted or obtained. Now we have the
following analogue of Morton's result. (We use again
$n=n(L)$ for the number of components, $g=g(L)$ for the genus and
$\chi=\chi(L)$ for the maximal Euler characteristic of $L$.)

\begin{prop}\label{i6,}
For $n\ge 4$ components, there are infinitely many (arborescent) 
canonically fibered links with any given monic admissible Alexander
polynomial.
\end{prop}

\proof We use the links of Theorem \reference{HHHH}. If $g>0$,
the unknotted component created by two claspings allows to apply
Stallings twists if we choose the claspings to be of opposite sign.
The linking number easily distinguishes infinitely many of the
resulting links, but they all have the same complements, so
hyperbolicity is preserved. For $g=0$ we can use Stallings twists for
the links in \eqref{uu} and for those of Gabai's type (C). (See the
proof of Theorem \reference{HHHH}.) His pretzel links of type (B)
are already infinitely many (and all have the same polynomial). \qed

We know in contrast (see the discussion at the end of this
section) that a generic monic Alexander
\em{knot} polynomial of degree 2 is realized by only finitely many
\em{canonical} fiber surfaces. So the combination of fibering and
canonicalness poses non-trivial restrictions on infinite families.
Assuming canonicalness and merely minimal genus property, the
scope of constructible infinite families widens.

\begin{prop}\label{pza}
For $n=1$ and $g>0$, or $n\ge 3$, any admissible Alexander link
polynomial $\Dl\ne 0$ is realized by infinitely many prime arborescent
$n$-component links with a canonical minimal genus surface and
$2\Md\Dl=1-\chi$.
\end{prop}

\proof For knots (and $\Dl\ne 1$) this can be shown by applying
the surgeries of the type \eqref{tgs} for all admissible $p,q,r$
at the parallel clasp * of the knots as in figure \ref{figgr},
constructed in the proof of Theorem \ref{thbd}. The distinction of
the resulting knots is a bit subtle, but since they are arborescent,
it can be done at least from \cite{BS}. For links of $\ge 3$
components and $g>0$, as in the proof of Theorem \reference{HHHH},
a parallel clasp is created by \eqref{clasping}, and the same surgery
applies. (For $n\ge 4$ the ``Stallings twist'' in proposition \ref{i6,}
would also apply, and the resulting links
are again much less sophisticatedly distinguished by linking
numbers.) The case $g=0$ and $n\ge 3$ is again easily recovered by
the pretzels.
\qed

For 2 components, however, some new idea is needed. The parallel
clasp disappears, and so far we cannot prove the claim, except
for special families of polynomials (it is also false if $g=0$).


Turning back to fiberedness, we do not know about extensions
of Morton's construction, explained in the beginning of  this
section, to obtain infinite families of links up to 3 components.
The infinite realizability is (even for general links or fiber
surfaces) not fully clear. As an application of our work we can
obtain at least the following additional examples. 

\begin{prop}\label{P73}
(1) For $n=3$ components and a monic admissible Conway polynomial
$\nb$ with $[\nb]_2=-1$, there exist infinitely many canonically
fibered links realizing $\nb$, which are connected sums of $2$ prime
arborescent factors. \\
(2) For knots ($n=1$), the same holds for polynomials $\nb$ with
a multiple zero. If $\nb=\nb_1^2$ for some $\nb_1\in\bZ[z]$, then
there exist infinitely many canonically fibered prime (arborescent)
knots realizing $\nb$.
\end{prop}

\proof For (1) take a prime fibered knot $K$ with $\nb_K=-z^{-2}
\nb(z)$, and build the connected sum with $(2,-2,2k)$-pretzel
links. Part (2) is an adaptation of the observation of Quach
\cite{Quach}. It suffices to consider the case $\nb=\nb_1^2$. If
$\nb_1\in\bZ[z]$ and $\nb_1^2\in\bZ[z^2]$, then $\nb_1\in\bZ[z^2]$
or $\nb_1\in z\bZ[z^2]$. Since $[\nb]_{z^0}=1$, former
alternative applies. Then w.l.o.g. $[\nb_1]_{z^0}=1$ up to
taking $-\nb_1$ for $\nb_1$.

So we can take a knot $K$
as in Theorem \reference{thbd} with $\nb_K=\nb_1$, and build the
connected sum $K\#!K$ at the parallel clasps in figure \ref{figgr}.
The canonical surface of the resulting diagram admits Stallings
twists at the spot of the connected sum. Since smoothing out a
crossing created by such Stallings twists gives a diagram of an
amphicheiral 2-component link $L$ (so that $\nb_L=0$), again
\eqref{sRs} shows that the twists preserve $\nb$. Also it is
easy to observe that the diagrams are still arborescent, so
infinitely many of the knots can be distinguished using \cite{BS}.
(There is again a much less sophisticated distinction argument,
which uses the leading term in the Alexander variable of the
skein polynomial.)
\qed

Since a fiber surface is connected, we must have $\chi\le 2-n$. For
$(n,\chi)=(2,0)$ we have only the Hopf links. For $(n,\chi)=(3,-1)$
and $\nb=+z^2$ we have again only 2 (composite) links, the connected
sum of two positive or two negative Hopf links (see part \ref{L51}
of Remark \reference{R51}). These observations are valid not only
for canonical, but also for general fiber surfaces, as is explained
in \cite{Kanenobu2}.

For $n=2$ and $\chi<0$, we can observe that the knots in
Morton's construction (see the proof of Theorem 4 in \cite{Morton})
likewise have unknotting number $1$, which allows to obtain
analogously to our case certain fibered 2-component links. It
seems some effort needed to extend Morton's JSJ decomposition
arguments and show that infinitely many of these links are
different. (Fibering and control of the Alexander polynomial
are again not difficult.) One would then have also (at least the
obvious connected sum) examples of 3 components for any polynomial.

We also do not know how to find for general (monic or not)
polynomials infinitely many (fibered or not) knots with certain
specific properties (like arborescent, prime, hyperbolic etc.).
For knots ($n=1$), part (2) of proposition \ref{P73} implies

\begin{corr}\label{cff}
In genus $g\ge 4$, then there exist infinitely many monic polynomials
realized by infinitely many canonically fibered prime knots. \qed
\end{corr}

To reformulate this more suitably, let for $d\ge 1$,
\[
\Phi_d\,:=\,\left\{\,\begin{tabular}{c}$\nb$ monic of degree $2d$,
realized by \\[1mm] infinitely many canonically fibered knots
\end{tabular}\,\right\}\,.
\]
Then we can understand $\Phi_d\subset \Gm_d:=\{\pm 1\}\times \bZ^{d-1}$.
We say that $\Phi_d$ is infinite if $d\ge 4$. Contrarily, $\Phi_1=\vn$,
and our aforementioned result in \cite{canon} shows that $\Phi_2$ is
finite. (We do not know about finiteness of $\Phi_3$.) So we see that,
expectedly, this result does not extend to $d\ge 4$, at least in
full strength. Nevertheless, for some $d$ still the inclusion
$\Phi_d\subset \Gm_d$ may be proper, or in fact so that $\Gm_d\sm
\Phi_d$ is infinite. The right sort of question to ask about what
polynomials are realized infinitely many times, seems to be something
like:

\begin{question}
Is $\Phi_d\subset \Gm_d$ contained
in the image of finitely many $d-1$-tuples of polynomials
\[
(f_1,\dots,f_{d-1})\in\bQ[x_1,\dots,x_k]^{\times d-1}\,,
\]
each $f_i$ of which maps $\bZ^k$ to $\bZ$, with $k\le d-2$?
\end{question}

There is a corresponding problem for links.
The question on the maximal $k$ needed also has some right.
The bound $d-2$ may be improvable, but obviously not below $1$
for $d=4,5$, and, with the origin of corollary \ref{cff} in mind,
expectably not below $d-4$ for $d\ge 6$.

\section{Large volume knots\label{S8}}

\subsection{Arborescent knots}

While so far we were concerned in estimating volume from above,
we give, using tangle surgeries, two constructions to obtain
knots of given polynomial and large volume. The case of links
is left out mainly for space (rather than methodological) reasons.
The first construction yields arborescent knots.

\begin{theorem}\label{theo2}
Given an Alexander knot polynomial $\Dl$ with $d=\Md\Dl$ and an integer
$g_s\ge \max(1,4d-1)$, there exist hyperbolic arborescent knots of
arbitrarily large volume with Alexander polynomial $\Dl$ and 4-genus
$g_s$.
\end{theorem}

Our result is motivated by similar work of Kalfagianni
\cite{Kalfagianni}, one of whose consequences (Corollary 1.1
therein) it improves. (At the end of this paper we will be able
to recover Kalfagianni's full result; our tools are, however,
somewhat different from hers.) A related result, that implies 
a certain part of the statement of Theorem \reference{theo2}, 
was obtained simultaneously by Silver and Whitten \cite{SWh}. 

%

%

\begin{lemma}\label{lmo}
The tangle surgeries \eqref{tgs} of Lemma \ref{lsurg} (for $k=1$)
alter $g_s$ most most by $\pm 2$.
\end{lemma}

\proof 
We like to examine the change of $g_s$ under the surgery.
We change first a crossing in the twist of $q$.
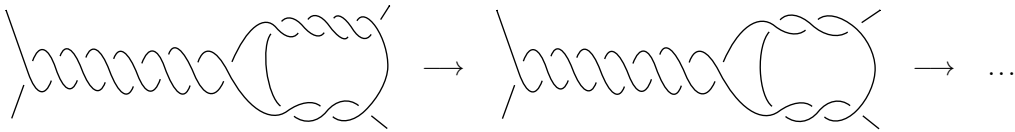
\[
\epsfsv{3cm}{t1-dlsurg1}\quad\lra\quad
\epsfsv{3cm}{t1-dlsurg2}\quad\lra\quad\dots
\]
Since $q+r=2$, applying concordance, we can cancel the remaining $q-2$
crossings with the crossings in the twist of $r$, and then remove the
(crossings in the) twist of $p$. Then by switching a crossing
we recreate the clasp before the tangle surgery.
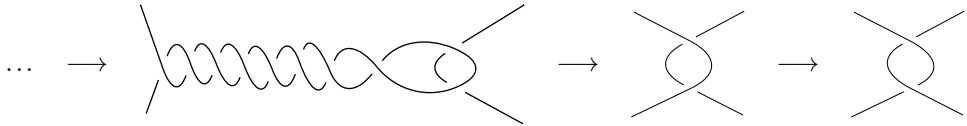
\[
\dots\quad\lra\quad\epsfsv{3cm}{t1-dlsurg3}\quad\lra\quad
\epsfsv{1.5cm}{t1-dlsurg4}\quad\lra\quad
\epsfsv{1.5cm}{t1-dlsurg5}
\]
Now $g_s$ changes by at most $\pm 1$ under a crossing change, so
it changes by at most $\pm 2$ under the tangle surgery. \qed

\proof[of Theorem \reference{theo2}]
In the following we choose integer triples $(p,q,r)$ with $p,q,-r>1$
odd, $r+q=2$ and $pq+pr+qr=-1$. We will assume that $p,q,r$ have
these properties throughout the proof. 

Choose from Theorem \reference{thbd} an arborescent knot $K$ with
$\Dl_K=\Dl$ and the arborescent diagram $\hat D$ constructed in the
proof. Following \cite{Adams} we call a crossing a \em{dealternator}
if it belongs to a set of crossings whose switch makes the diagram
alternating. This set is determined up to taking the complement.
Since we constructed $\hat D$ to have at most $4d-2$ twist
equivalence classes, we can choose (possibly taking the complement)
the number $d$ of twists in $\hat D$ consisting of dealternators to be 
\[
t\le 2d-1\,.
\]

Now we can turn $\hat D$ into an arborescent diagram $\hat D_0$ of
$K$, so that each of the $d$ twist equivalence classes of
dealternators in $\hat D$ becomes a single (dealternator) crossing
in $\hat D_0$. Fix in $\hat D_0$ the set of $d$ dealternators so
obtained. Create (by a Reidemeister II move) a trivial parallel
clasp near each dealternator, obtaining a diagram $D_0'$ of
$K$ with dealternators occurring in $d$ parallel clasps.

Now let $T_{p,q,r}$ be the tangle described in the proof of Lemma
\ref{lsurg} for $k=1$, and $T_{-p,-q,-r}$ its mirror image.
(So by the index shift $p$ means now what was $p+1$ in that proof.)
Let $D_0=D_0(p,q,r)$ be the result of substituting $T_{p,q,r}$
for each positive dealternator clasp tangle, and $T_{-p,-q,-r}$
for each negative dealternator clasp tangle in $D_0'$. Let
$K_{p,q,r}$ be the knot $D_0$ represents. Then $D_0$ has all
its dealternators in twists in the substituted tangles. When
now the length of the twists in $T_{p,q,r}$ grows, Thurston's 
hyperbolic surgery theorem shows that $\vol(K_{p,q,r})$ converges
(from below) to the volume of a certain link $T_\infty$. This limit
link is the same as when $r$ has opposite sign, but then we have
prime alternating diagrams. So $T_\infty$ is an augmented
alternating link (as in \cite{Brittenham,Lackenby}). Then in order
to obtain large volume we apply Adams' result on the volume of
augmented alternating links (see \cite{Brittenham,Lackenby}), and
so it is enough to increase the number of tangles whose twist lengths
we can augment unboundedly.



Simultaneously we want to carry out our construction so as to
obtain large $g_s$. With $p,q,r$ given, we applied the tangle 
surgeries of Lemma \ref{lsurg} (for $k=1$) at each clasp of
dealternators in $D'_0$ and obtained a diagram $D_0=D_0(p,q,r)$.
By Lemma \reference{lmo} we have 
\begin{eqn}\label{est1}
|g_s(D_0)-g_s(K)|\,\le\,2t\,\le\,4d-2\,.
\end{eqn}
Since $u(K)=1$, we have $g_s(K)\le 1$, so $g_s(D_0)\le 4d-1$.

We consider the pretzel knots $P(p,q,r)$, which have $\Dl=1$.
By the main theorem in \S 1 of \cite{Rudolph}, these 
pretzel knots are quasipositive, and by Proposition 5.3
of \cite{Rudolph} have slice genus $1$.

Let now $D=D(l,p,q,r)$ be the diagram obtained by taking connected 
sum of $D_0$ with $l$ copies of the $(p,q,r)$-pretzel diagram. 
(Note that now $p,q,r$ enter into the construction of $D(l,p,q,r)$
in a second different way.) Because $P(p,q,r)$ is quasipositive of
$4$-genus one, we have by the Bennequin-Rudolph inequality
(see \cite{Rudolph2})
\begin{eqn}\label{est2}
g_s(D(l,p_l,q_l,r_l))\to\infty
\end{eqn}
when $l\to\infty$, for any sequence $(p_l,q_l,r_l)$ of
triples $(p,q,r)$ of the above type. Moreover, the
numbers \eqref{est2}, when taken over all $l\ge 0$,
realize all integers $g_s\ge 4d-1$, again regardless of the choice of
$(p_l,q_l,r_l)$.

We apply now the moves \eqref{tsurg}.
Choose the connected sum in $D$ so that the creation of two parallel
clasps in the first move in \eqref{tsurg} gives a prime diagram $D'$.
The second move is a tangle surgery, which preserves $\Dl$
and can be performed for any triple $(p,q,r)$. (In \eqref{tsurg}
we show the operation for the simplest triple, which after the
shift of $p$ is now $(7,5,-3)$.) Call the resulting diagram
$D''=D''(l,p,q,r)$, and $K''=K''(l,p,q,r)$ the knot it represents.
Since this surgery is a concordance, we have 
\begin{eqn}\label{est3}
g_s(D'')=g_s(D)\,.
\end{eqn}

So from \eqref{est2} and \eqref{est3} we have then
\[
g_s(D''(l,p_l,q_l,r_l))\to \infty\,,
\]
when $l\to\infty$ and $(p_l,q_l,r_l)$ is an arbitrary
sequence of tuples $(p,q,r)$. Moreover, all numbers above
or equal to $4d-1$ are realized as $4$-genera. Now $D''$ has all its
dealternators occurring in twists whose length can be augmented
arbitrarily, preserving $\Dl$. So if for each $l$ we choose
$-r_l$ (and hence $q_l,p_l$) large enough, we obtain hyperbolic
knots $K_l=K''(l,p_l,q_l,r_l)$ of large volume from the results
of Thurston and Adams. 

In order to obtain infinitely many knots of fixed 4-genus
take in the construction of $D(l,p,q,r)$ connected sum with
$(p,q,r)$-pretzel diagrams and mirror images thereof (with
reverse orientation). The volume will distinguish infinitely
many of the knots $K_l$.

To verify that $K_l$ is arborescent, use that we chose the initial
diagram $\hat D$ of $K$ to be arborescent. Taking iterated 
connected sum with the $(p_l,q_l,r_l)$-pretzel knots and adding clasps 
can be done so as to preserve arborescency of the diagram. The same
observation applies to the tangle surgeries. \qed

Using the upper bound in Theorem \ref{LAT}, we have a result on 
growing twist numbers. 

\begin{corr}
Any possible Alexander polynomial is realized by arborescent knots
$K_l$ with twist number $t(K_l)\to\infty$. \qed
\end{corr}

\begin{rem}\label{rSW}
Our construction can be easily adapted to preserve the Alexander module.
Choose a prime $s$ such that all (finitely many up to units) divisors
of $\Dl$ in $\bZ[t^{\pm1}]$ (including $\Dl$ and $1$) remain distinct
(up to units) when coefficients are reduced $\bmod s$.  Then choose
$p,q,r$ so that $p+1,q,r\equiv 1\,(2s)$, by choosing (for $k=1$)
$n$ in \eqref{xxx} divisible by $s$. Observe that changing any of
$p,q,r$ by (multiples of) $2s$ preserves a (properly chosen) Seifert
matrix $\bmod s$, and the Seifert matrix determines the Alexander
module. Since our arguments incorporate concordance, we can recover
most of the properties obtained by Silver and Whitten \cite{SWh},
except of course the knot group homomorphism. 
\end{rem}



\subsection{Free genus}

Our final result combines all the methods introduced previously
to obtain an extension of a theorem of Brittenham \cite{Brittenham2}.
He constructed knots of free genus one and arbitrary large volume.
We state a similar property for free genus greater than one.

\begin{theorem}\label{thfg}
Let $\Dl$ be an admissible Alexander knot polynomial of degree
$d\ge 2$. Then there exist hyperbolic knots $K_n$ of arbitrarily
large volume with free genus $g_f(K_n)=d$ and $\Dl(K_n)=\Dl$.
\end{theorem}

\begin{rem}
As to extensions and modifications of this statement,
the following can be said:
\begin{mylist}{\arabic}
\myitem Our construction does not apply for free genus one.
 The Alexander polynomial is not of particular interest on
 genus one knots, so its control in Brittenham's (or some
 similar) construction seems only of minor use, and we will not
 dwell upon this here.
\myitem A justified question is whether for monic polynomial we
 can actually find fibered knots. We expect that it is possible,
 but the effort of proof would grow further, too much for the
 intention and length of this paper.
\myitem Another suggestive question, whether one can replace
 free by canonical genus, is to be answered negatively.
 Brittenham had shown \cite{Brittenham} that canonical genus
 bounds the volume (see also \cite{gen2}).
\myitem The knots we obtain are unlikely arborescent or of
 unknotting number one, but still have slice genus at most one
 if $g_f>2$.
\myitem The case of links is, like the explanation at the
 beginning of this section, analogous to treat (with similar
 mild constraints), but also left out for space reasons.
\end{mylist}
\end{rem}

\proof Let first $g_f>2$. For a given number $k$ we consider
the link $L_k=K\cup U_1\cup\dots\cup U_k$ given by replacing
the diagram of the knot $K$ from theorem \ref{thbd} along the (more 
tightly) dashed line $\gm$ below as follows:\\[1mm]
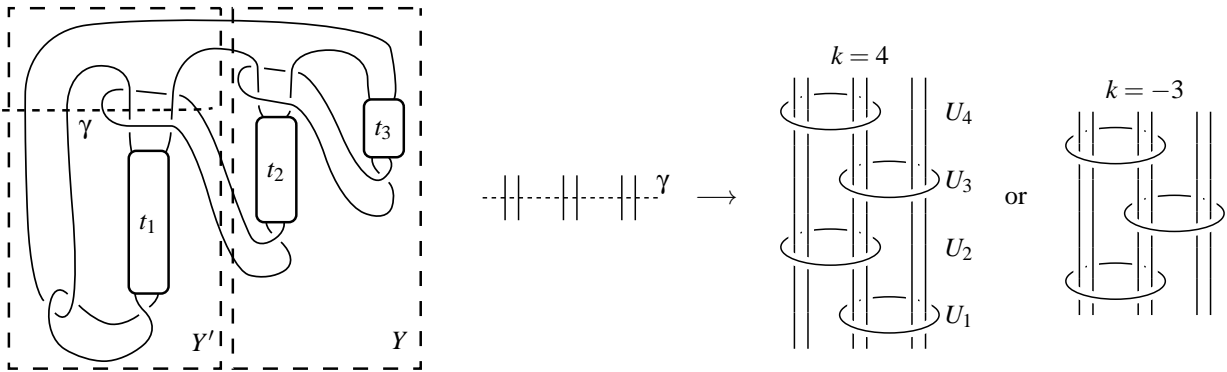
\[
\vcbox{\epsfdiag{%
  \kern3mm\raisebox{2mm}{\epsfs{5cm}{t1-alrealiz}}}{1cm}{
    \pictranslate{0.3 0.2}{
       \picfilledovalbox{1.65 1.85}{0.53 1.9}{0.1}{$t_1$}
       \picfilledovalbox{3.35 2.55}{0.53 1.4}{0.1}{$t_2$}
       \picfilledovalbox{4.78 3.10}{0.53 0.77}{0.1}{$t_3$}
       {\piclinedash{0.1}{0.00}
        \opencurvepath{-0.3 3.34}{1.3 3.34}{2.4 3.35}
	  {2.5 3.37}{} 
       }
       \picputtext{2.38 0.27}{$Y'$}
       \picputtext{5. 0.27}{$Y$}
       \picputtext{0.80 3.13}{$\gm$}
       {\piclinedash{0.2}{0.1}
       \picbox{1.2 2.3}{2.82 4.8}{}
       \picbox{4.02 2.3}{2.5 4.8}{}
       }
    }
  }
}
\kern1cm
\diag{6mm}{4}{1}{
  {\piclinedash{0.1}{0.00}
   \picline{0 0.5}{4 0.5}
  }
  \picmultigraphics{3}{1.3 0}{
    \picline{0.5 0}{0.5 1}
    \picline{0.8 0}{0.8 1}
  }
  \picputtext{4. 0.87}{$\gm$}
}\quad\lra\quad
\def\ring{
      \picmultigraphics{3}{1.3 0}{
	\picline{-1.45 0}{-1.45 0.8}
	\picline{-1.15 0}{-1.15 0.8}
      }
      \picmultiellipse{-8 1 -1 0}{0.65 0.75}{1.1 0.4}
      \picmultigraphics{3}{1.3 0}{
	\picmultiline{-8 1 -1 0}{-1.45 1.5}{-1.45 0.8}
	\picmultiline{-8 1 -1 0}{-1.15 1.5}{-1.15 0.8}
      }
}
\diag{6mm}{4.6}{6.7}{
  \pictranslate{2 0}{
    \picmultigraphics{2}{0 3}{ \ring }
    \picscale{-1 1}{\pictranslate{0 1.5}{
      \picmultigraphics{2}{0 3}{ \ring }}}
    \picputtext{0 6.5}{$k=4$}
    \picputtext{2.2 .75}{$U_1$}
    \picputtext{2.2 3.75}{$U_3$}
    \picputtext{2.2 2.25}{$U_2$}
    \picputtext{2.2 5.25}{$U_4$}
  }
}\quad\mbox{or}\quad
\diag{6mm}{4}{5.2}{
  \pictranslate{2 0}{
    \picmultigraphics{2}{0 3}{ 
      \picscale{-1 1}{ \ring }
    }
    \pictranslate{0 1.5}{ \ring }
    \picputtext{0 5}{$k=-3$}
  }
}
\]
(We extend this to $k<0$ by placing the circles $U_i$ the other 
way, as shown.) Choosing $a,b$ sufficiently large, we construct 
the knots $K_{k,a,b}$ from $L_k$ for $4\mid k$ by doing 
\[
\left\{
\begin{array}{c}
a \\ b \\ -a \\ -b
\end{array}\right\}\mbox{ twists at $U_k$ for $k\equiv $ }
\left\{
\begin{array}{c}
1 \\ 2 \\ 3 \\ 0
\end{array}
\right\}\,\bmod 4,
\mbox{\ in the following way:}
\diag{8mm}{3}{1.5}{
  \pictranslate{1.5 0}{
    \picmultigraphics{2}{1.3 0}{
      \picline{-.8 0}{-.8 1.5 2 :}
      \picline{-.5 0}{-.5 1.5 2 :}
    }
    \picmultiellipse{-8 1 -1 0}{0. 0.75}{1.1 0.4}
    \picmultigraphics{2}{1.3 0}{
      \picmultiline{-5 1 -1 0}{-.8 3 4 :}{-.8 1.5}
      \picmultiline{-5 1 -1 0}{-.5 3 4 :}{-.5 1.5}
    }
    \picputtext{1.2 1.3}{$+1$}
  }
}
\quad\lra\quad
\diag{8mm}{2.5}{3}{
  \pictranslate{1.25 0}{
    \picmultigraphics{2}{0 1.5}{ 
      \picmultigraphics[S]{2}{-1 1}{
	\picmultigraphics{2}{0.4 0}{
	  \picmulticurve{-5 1 -1 0}{0.5 0}{0.5 0.5}{-0.9 1.0}{-0.9 1.5}
	}
      }
    }
  }
} \,.
\]

Here a few annotations seem proper. (i) The twists along $U_k$ 
are called in the common cut-paste-language surgeries. However,
we avoid this term here in order not to confuse with the tangle
surgeries (which will just reenter). The ``twists'' may, in turn, 
conflict with definition \ref{df22}, but they can be regarded here
as an extension of the previous concept, and so seem the more
convenient term. (ii) Twisting along $U_i$ adds also a full
twist (now in a sense directly related to definition \ref{df22})
into the bands. However, these twists cancel each other when
twisting at $U_i$ is performed in the prescribed way, so we can
ignore them.

It is easy to see now that $K_{k,a,b}$ has the same Alexander 
polynomial as $K$, since the Seifert matrix is not altered by 
the twisting at $U_i$. Similarly, the twisted Seifert surface
is still free. By thickening the surface into a bicolar, we
see that the twisting at $U_i$ accounts only in braiding the 
various 1-handles, and this braiding can be undone by sliding 
the handles properly, as for the braidzel surfaces 
\cite{Rudolph,Nakamura3}. 

With this we focus on hyperbolicity. By Thurston and Adams
again it suffices to show that $L_k$ are hyperbolic for large 
$|k|$. (We need in fact here only $k>0$ and $4\mid k$, but we 
will soon see why it is good to have the other $k$ around, 
too.) We use the tangle decomposition $Y\cup Y'$ of $K$, which 
carries over with modifications to $L_k$. (In order not to
overwork notation, we denote $Y,Y'$ the same way in all 
links, each time specifying the link.) First we use tangle 
surgery to remove the dependence of $Y'$ in $L_k$ on the number 
$t_1$ of full twists. The surgery allows us to replace the 
lower part of $Y'$ as follows:
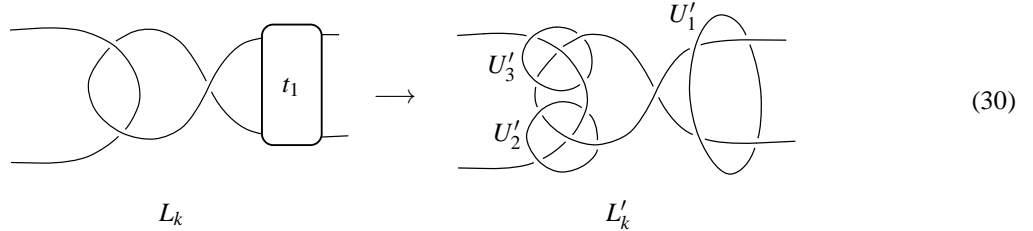
\begin{eqn}\label{acf}
\begin{array}{ccc}
\\[1mm]
\vcbox{
\epsfdiag{\epsfs{3cm}{t1-alrealiz3}}{1.5cm}{
  \piclinewidth{60}
  \picfilledovalbox{2.5 0.7}{0.53 1.07}{0.1}{$t_1$}
}}
& \lra &
\vcbox{
\epsfdiag{\epsfs{3cm}{t1-alrealiz3a}}{1.5cm}{
  \piclinewidth{40}
  \picputtext{2.01 1.45}{$U_1'$}
  \picputtext{0.45 0.4}{$U_2'$}
  \picputtext{0.40 1.0}{$U_3'$}
}}
\\
\ry{1.6em}L_k & & L_k'
\end{array}
\end{eqn}
The meaning is that we can have a free surface and a desired 
Alexander polynomial by applying a proper, but arbitrarily
augmentable, number of twists at the circles we added. Now 
$U_1'$ is in fact parallel to $U_1$ for $k=1$. So we can, and
for hyperbolicity must, omit $U_1'$ then. This can be done 
with the understanding that we perform at $U_1$ the additional
twists we would have needed to perform at $U_1'$.

The effect of the surgery is now that the link $L_k'$,
whose hyperbolicity it suffices to show, has a tangle $Y'$ 
which does no longer depend on $t_1$, but only on $k$. 

\begin{lemma}
The links $L=L_k'$ are prime.
\end{lemma}

\proof The (only, but then indeed so because of $\Dl$) knotted 
component $K$ of $L_k'$ is prime; e.g. it has unknotting number 
one. Thus if $L$ is composite, there is a composite (possibly 
split) 2-component sublink $L'=K\cup O$ of $L$. Now, for such
sublinks, the tangle $Y'$ reduces only to finitely many cases;
in fact 3 are enough to test (using that $U_1$ and $U_1'$
in \eqref{acf} are parallel, and $U_2'$ and $U_3'$ are
flype-equivalent). These 4 tangles $Y'$ can be checked to 
be prime by \cite{KL}, and since the same can be done for $Y$ 
(despite of its dependence on $t_2,t_3,\dots$), we have that 
$L'$ is prime, a contradiction. Thus $L$ is prime.
\qed

\begin{lemma}
The links $L=L_k'$ are atoroidal.
\end{lemma}

\proof We first prove for $|k|\le 3$. The main point here is
to remove the dependence of $Y$ on $t_2,t_3,\dots$.

The $t_2$ twists can be easily removed by tangle equivalence. The 
argument that eliminates $t_3,t_4,\dots$ consists in a repetition
of our work in applying Oertel's and Wu's results, so we just 
recapitulate the main points. Now $Y_1=Y$ and $Y_2=Y'$, and we
have $L=Y\cup Y'$, with $Y$ 
being a Montesinos tangle of length 2 for $g_f(K)=3$, or an 
arborescent tangle subjectable to Wu's result for $g_f(K)>3$.
Assume $T$ is an essential torus of $L$. Then again $T\cap X(Y)$
is empty, all of $T$, or an annulus $A$. 

If $A$ exists, then by Sublemma \reference{T1k} and the argument
after it, $A$ is $\pa$-parallel to $C$, so can be moved out. If
$T\subset X(Y)$, we have a contradiction to Wu for $g_f>3$, or
by gluing $Y$ and $A$ to itself and Oertel's result if $g_f=3$.

If $T\subset X(Y')$, then $T$ is essential in $E(L)$ even after
modifying $Y$, as long as $Y$ is prime and has a closed component.
Since for $|k|\le 3$, we have only finitely many $Y'$, we can
easily find a proper prime tangle $Y$ and check the hyperbolicity
of the handful of links $L=Y\cup Y'$ by SnapPea to see the
contradiction to the existence of $T$. With this argument 
the atoroidality is proved for $|k|\le 3$. 

Now let $|k|\ge 4$. We use induction on $|k|$ (where
the cases $4\nmid k$ enter). Assume $T$ 
is again an essential torus of $L'_k$. By induction, $T$ is 
inessential in $L'_k\sm U_1$. (Here the use of $L'_k$ also for 
$k<0$ pays off.) So $T$ contains in one of its complementary 
regions either only $U_1$, or $U_1$ and exactly one other 
component $V$, to which it becomes $\pa$-parallel after 
removing $U_1$. In particular in latter case $T$ must have 
the knot type of $V$. Applying the same argument to $L'_k\sm 
U_{|k|}$ shows then that $T$ must contain exactly $U_1$ and 
$U_{|k|}$ in one of its regions $R$, be unknotted, and have 
them as cores of the solid torus $R=:\ir T$. 

But now, if one removes $U_2$ and $U_3$ from $L$ (here the
assumption $|k|\ge 4$ enters), then again $T$ must become
inessential. However, $\wt L=L\sm U_{2,3}$ is non-split by the 
previous lemma, and the exterior of $T$ in $E(\wt L)$ contains 
the knotted component $K$. Then $T$ cannot compress or be 
$\pa$-parallel in its exterior, but the same applies to its 
interior either, a contradiction. With this the lemma is 
proved. \qed

\begin{lemma}
The links $L=L_k'$ are not Seifert fibered.
\end{lemma}

\proof Again components of Seifert fibered links are (possibly
trivial) torus knots, and for our links we have a knotted component 
of unknotting number one. It must be then a trefoil, but then 
we are in the situation $g_f=1$, which we chose not to consider. 
\qed

Now we have shown the theorem for $g_f>2$. Our procedure does not
work, though, for $g_f=2$ (exactly the same way; for example, then
$Y$ is no longer prime). In that case, we realize $V_2$ of \eqref{57}
as a Seifert matrix in the way shown in the diagram (a)
of \eqref{gf2}. Here we 
took the example with $a_1=a_2=2$. The Conway polynomial is $\nb=
1-a_1z^2+a_2z^4=1-2z^2+2z^4$. In general the half-twists
at * are $2a_1-1$, and those at ** are $2a_2+1$. (Again
$-1$ half-twist is a crossing of negative skein sign.)

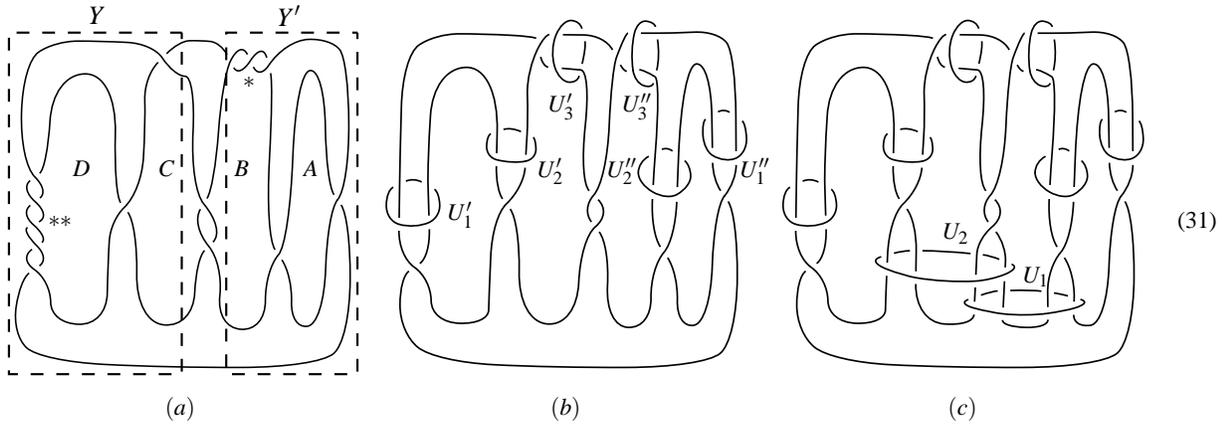
\begin{eqn}\label{gf2}
\small
\begin{array}{c@{\quad\es}c@{\quad\es}c}
\epsfdiag{\epsfs{4.5cm}{t1-alrealiz4}}{1.42cm}{
  \picputtext{2.22 2.68}{$*$}
  \picputtext{0.45 1.37}{$**$}
  \picputtext{0.65 1.87}{$D$}
  \picputtext{1.45 1.87}{$C$}
  \picputtext{2.15 1.87}{$B$}
  \picputtext{2.8 1.87}{$A$}
  \scalelines
  \picscale{0.6 d}{\piclinedash{0.2}{0.1}
   \picbox{1.3 2.57}{2.70 5.34}{}
   \picbox{4.37 2.57}{2.05 5.34}{}
  }
  {\normalsize
  \picputtext{0.8 3.3}{$Y$}
  \picputtext{2.6 3.3}{$Y'$}
  }
} &
\epsfdiag{\epsfs{4.8cm}{t1-alrealiz4a_0}}{1.44cm}{
  \picputtext{2.32 2.45}{$U_3''$}
  \picputtext{1.62 2.45}{$U_3'$}
  \picputtext{0.69 1.45}{$U_1'$}
  \picputtext{3.41 1.85}{$U_1''$}
  \picputtext{1.52 1.85}{$U_2'$}
  \picputtext{2.18 1.85}{$U_2''$}
} &
\epsfdiag{\epsfs{4.8cm}{t1-alrealiz4a}}{1.44cm}{
  \picputtext{1.6 1.26}{$U_2$}
  \picputtext{2.34 0.86}{$U_1$}
} \\
\ry{1.6em}(a) & (b) & (c)
\end{array}
\end{eqn}

The rows/columns of $V_2$ correspond to curves that go in 
positive direction along the regions $A,B,C,D$. The curves
for $A$ and $B$, resp. $C$ and $D$, intersect once on the
lower Seifert circle; otherwise curves do not intersect.

Now observe that again we can apply a surgery in $Y$ and $Y'$ 
(where in lemma \ref{lsurg}, we have $k=1$ for $Y'$ and $k=a_2\ne 0$ 
for $Y$). It allows to arbitrarily augment the number of twists,
keeping $\Dl$ and the surface canonical. This has the effect
of eliminating the dependence on $\Dl$ (i.e. on $a_{1,2}$) of 
the link, whose hyperbolicity it is enough to show; see (b) 
in \eqref{gf2}. Denote the triples of circles occurring for 
the surgery in $Y$ by $U_i'$, and let those for $Y'$ be $U_i''$. 

Finally, we must add the circles $U_i$ around pairs of bands.
This is done as shown for $k=2$ in part (c) of \eqref{gf2}.
Since the links $L_k$ we obtain depend only 
on $k$, we can use the same type of inductive argument to
show atoroidality, checking the initial links by SnapPea.

We use then twisting at the $U_i$ again for $4\mid k$ in the previously
specified way. It may be worth remarking that, to see the preservance
of $\Dl$, the twists along $U_i'$ and $U_i''$, resulting from 
the tangle surgeries, must be performed before those at $U_i$.
The $U_i$ enter into the tangle the surgeries are performed at. 
Inspite of this, the resulting modifications are independent 
from each other, so no conflict arises.

To exclude a Seifert fibration for $L_k$, note that if the
not obviously unknotted component $K'$ is indeed knotted,
none of the Burde-Murasugi links has such a component (even
if a torus knot), and more than two unknotted ones. If $K'$
is unknotted, the Seifert fibration for $L_k$ is excluded
using linking numbers. A look at the Burde-Murasugi list shows
that there is no link with all linking numbers zero, except
the trivial link (unlink). This is excluded by looking at a
proper sublink of $L_k$.
\qed

\begin{rem}\label{rKf}
Observe that the twisting at the
components $U_i$ corresponds in an obvious way to a (power of
the) commutator $[\sg_1^a,\sg_2^b]=\sg_1^a\sg_2^b\sg_1^{-a}
\sg_2^{-b}$ in the 3-strand braid group $B_3$. Using higher
order commutators (and leaving out the tangle surgeries), one
can preserve, additionally to $\Dl$, Vassiliev invariants of
given degree. Then from the argument for $K'$ being unknotted,
one easily recovers the main result
of Kalfagianni \cite{Kalfagianni}: given $n>0$, there exist
hyperbolic knots $K_n$ of arbitrary large volume with $\Dl=1$
and trivial Vassiliev invariants of degree $\le n$. (In our
construction also $g_f(K_n)\le 2$.)
\end{rem}

It appears a good challenge, and may be a future investigation,
to extend this result by showing that $K_n$ can be chosen to
be $n$-similar (i.e. with Vassiliev invariants of degree $\le n$
coinciding) and with the same Alexander polynomial as any given
knot $K$. Kalfagianni's theorem is the statement for $K$ being
the unknot. Certainly our method offers more than just this special
case. For example, without keeping genus minimality of the surface,
one could easily find $K_n$ when $K$ is any Montesinos knot.

\begin{table}[ptb]
\captionwidth0.8\vsize\relax
\newpage
\vbox to \textheight{\vfil
\tabcolsep5.8pt
\rottab{%
\def\hh{\\[0.8mm]\hline[1.5mm]}%
\def\PBS#1{\let\temp=\\\centering #1\let\\=\temp}%
\hbox to \textheight{\hss
\def\rr{\rule[-4mm]{0mm}{0mm}}%
\def\eqq{{\small$2\Md\Dl=1-\chi_c$}}%
\def\yes{{\bf yes}\ }%
\def\hyp{{\bf hyp.}\ }%
\def\arbor{{\bf arbor.}}%
\begin{mytab}{|p{20mm}|@{\kern1mm}|
  p{28mm}|p{28mm}|p{28mm}|p{43mm}|p{36mm}|}
  { & & & & & }%
\hline [1mm]%
  \PBS{\mbox{}\\\# comps} &
    \PBS{arbitrary $\Dl\ne 0$ \\ \eqq \\ one link}
    &
    \PBS{arbitrary $\Dl\ne 0$\\ \eqq \\ $\infty$ many
    \rr}
    & 
    \PBS{monic $\Dl$\\ one canon.\\fibered link}
    &
    \PBS{monic $\Dl$\\ $\infty$ many canon.\\fibered links}
    &
    \PBS{monic $\Dl$\\ $\infty$ many\\fibered links}%
\\
\hline[-2mm]
& & & & & \\
\hline [2mm]%
\PBS{\mbox{}\\1} & \PBS{\yes (\arbor; Theorem \ref{thbd})\\
      \hyp for $g>0$ \\ (Remark \reference{zas})} &
    \PBS{\yes (\arbor; propos.\\ \reference{pza}) for $g>0$\\
      (no for $g=0$)} &
    \PBS{\yes (\arbor), \hyp\\ except unknot or trefoil\\
      (Theorem \reference{thbd})} &
   \PBS{no for $g\le 1$ and almost all\\ $\Dl$ in $g=2$
      \cite{canon};\\ \yes for $\nb$ with double\\
      zero (propos. \reference{P73}); \\
      unknown in general for $g\ge 3$\rr
   } &
   \PBS{no for $g\le 1$; yes\\ for $g\ge 2$ \cite{Morton}} 
   \\
\hline[-2mm]
& & & & & \\
\hline [2mm]%
\PBS{\mbox{}\\2} & \PBS{\yes (\arbor)\\\hyp for $g>0$ (theorem
     \reference{th2c})} &
    \PBS{unknown; no for $g=0$} & 
    \PBS{\yes (\arbor)\\\hyp for $g>0$ (Theorem \reference{th2c})} &
    \PBS{no for $g=0$ \\and almost
      all\\ $\Dl$ in $g=1$ \cite{canon};\\ else unknown} &
    \PBS{no for $g=0$; unknown,\\ likely yes (modif. of\\
     Morton; see \S\reference{S72}) if $g>0$\rr}
   \\
\hline[-2mm]
& & & & & \\
\hline [2mm]%
\PBS{\mbox{}\\3} & \PBS{\yes (part \ref{i4,}\\
    of Remark \ref{R51},\\ proposition
    \ref{pza};\\ \hyp \arbor)} & \PBS{\yes (propos.
    \reference{pza}; \arbor)} &
    \PBS{\yes (Theorem \reference{HHHH}; \\\hyp \arbor)
    if $\nb\ne +z^2$; \\ only compos. exist if $\nb=+z^2$} &
    \PBS{\yes if $[\nb]_2=-1$\\ (propos. \reference{P73};
    compos. links);\\  no if $\nb=z^2$, else unknown\rr} &
    \PBS{no if $g=0$, $\nb=z^2$ (see\\rem. in \cite{Kanenobu2});
      \yes if\\ $[\nb]_2=-1$ (compos.\\ links); else unknown\rr} \\
\hline[-2mm]
& & & & & \\
\hline [2mm]%
\PBS{\raisebox{-0.6em}{$\ge 4$}} &
  \PBS{\yes (\hyp \arbor)} & \PBS{\yes (\hyp \arbor)\rr} &
  \PBS{\yes (Theorem \ref{HHHH}; \hyp \arbor)} & \PBS{\yes (prop.\\
    \ref{i6,}; \hyp \arbor)\rr} &
    \PBS{\yes (\arbor)} 
\hline[-2mm]
& & & & & \\
\hline%
\end{mytab}%
\hss}%
}{The realizability status of given Alexander polynomials by
given number of given type of knots or links. The boldfaced
entries refer to the contribution of this paper.\protect\label{tab2}}
\vss}
\newpage
\end{table}

\section{
Questions and problems}

We mentioned already, for example in sections \ref{S72} and \ref{S8},
several problems, that may be the topic of future research.
We conclude with one other group of further-going questions,
concerning special knots realizing Alexander polynomials.

After we were able to incorporate arborescency into most of
our constructions, it makes sense to ask in how far one can
further restrict the type of knots.

\begin{question}
Are arbitrary Alexander polynomials realizable by Montesinos 
knots (perhaps), or even general pretzel knots (unlikely)?
\end{question}

The following argument shows that at least among pretzel knots 
restrictions on the Alexander polynomial may apply.

\begin{prop}
There exist Alexander polynomials not realizable by any generalized
pretzel knot $(a_1,\dots,a_{2n+1})$ with $a_k$ odd, for any $n$.
\end{prop}

\proof If we use equivalently $\nb$, then a direct skein argument shows
that all coefficients $\nb_j=[\nb]_{z^j}$ for even $j$, are polynomials
in $a_1,\dots,a_{2n+1}$ of degree at most $j$. (One can also 
argue with the work in \cite{bseq} and the well-known fact that
$\nb_j$ is a Vassiliev invariant of degree at most $j$.) Also, these
polynomials are at most linear in any $a_k$. Furthermore, they
are symmetric in all $a_k$, since permuting $a_k$ accounts for
mutations, that preserve $\nb$. So $\nb_j$ is a linear combination
of elementary symmetric polynomials $\sg_i$ in $a_k$ for $i\le j$.
Then one also finds that $\sg_j$ indeed
occurs in this linear combination, and only $\sg_i$ for even
$i$ occur. (Latter property is due to the fact that $\nb$ is
invariant under taking the mirror image.) So, up to linear
transformations, it is enough to see that some integer tuples
$(\sg_2,\sg_4,\dots, \sg_j)$, even for $\sg_i$ satisfying certain
congruences, cannot be realized as values of elementary symmetric
polynomials of any odd number of odd integers $a_k$. But $\sg_i$
occur as coefficients of the polynomial
\[
X(x)=(x-a_1)(x-a_2)\dots(x-a_{2n+1})\,,
\]
and it is known that the coefficients of polynomials with
real roots satisfy certain inequalities; they are $\log$-concave
(see Theorem 53 in \cite{HLP}).
So for example any triple $(\sg_2,\sg_4,\sg_6)$ with
$0<\sg_4<\sg_2<\sg_6$ will not occur. \qed

Another question addresses an important
point as to how a volume estimate can be strengthened.

\begin{question}
Is there a global constant $C$, such that all Alexander 
polynomials are realized by hyperbolic knots of volume $\le C$?
\end{question}

One can pose the analogous questions also for links.

%
%

\section{Result summary}

Table \ref{tab2} summarizes the state of knowledge about
realizing (monic) Alexander polynomials by links with a
canonical minimal genus (or fiber) surface, depending on
the number of components, the Alexander polynomial and
whether one or infinitely many such links are sought.

\noindent{\bf Acknowledgements.}
I would like to thank to Mikami Hirasawa, Efstratia Kalfagianni,
Taizo Kanenobu, Kunio Murasugi, Takuji Nakamura, Makoto Sakuma,
Dan Silver and Ying-Qing Wu for some helpful remarks, discussions,
and references. This work was partly carried out under Japan
Society for the Promotion of Science (JSPS) Postdoc grant P04300
at the Graduate School of Mathematical Sciences, University of Tokyo.
I also wish to thank to my host Prof. T.~Kohno for his support.

	
\end{document}